\documentclass{eccm-ecfd}
\usepackage{graphicx}
\usepackage{amsmath}
\usepackage{amsfonts}
\usepackage{amssymb}
\usepackage{latexsym}
\usepackage{amssymb}
\usepackage{amsbsy}
\usepackage{amstext}
\usepackage{amsxtra}
\usepackage{amscd}
\usepackage{amsopn}
\usepackage{lmodern}
\usepackage{ulem}
\usepackage{appendix} 
\usepackage{upgreek}
\allowdisplaybreaks
\usepackage{graphicx}
\usepackage[table,xcdraw]{xcolor}
\usepackage{sectsty}
\usepackage{multirow}
\usepackage[T1]{fontenc} 
\usepackage[nodisplayskipstretch]{setspace}
\usepackage{tabu}
\usepackage{booktabs}
\usepackage{array}
\usepackage{microtype} 
\usepackage[english]{babel}
\usepackage{titlesec}
\usepackage{eufrak}
\usepackage{nomencl}
\makenomenclature
\usepackage{hyperref}
\usepackage{mathtools}
\usepackage{placeins}
\usepackage{fancybox,framed}
\usepackage{listings}
\usepackage{color,soul}
\definecolor{light-gray}{gray}{0.8}
\sethlcolor{light-gray}

\title{On the derivatives of curvature of framed space curve and their time-updating scheme: Extended version with MATLAB code}

\author{MAYANK CHADHA$^{1,}$\footnote{Corresponding author.} AND MICHAEL D. TODD$^{2}$}

\heading{Curvature and its derivatives: Mayank Chadha and Michael D. Todd}

\address{$^{1}$ University of California, San Diego\\
9500 Gilman Drive, La Jolla CA 92092-0085\\
machadha@eng.ucsd.edu
\and
$^{2}$ University of California, San Diego\\
9500 Gilman Drive, La Jolla CA 92092-0085\\
mdtodd@eng.ucsd.edu\and
This is the extended version of paper submitted to \textit{Applied Mathematics Letters} Journal, Elsevier with the title ``\textit{On the derivatives of curvature of framed space curve and their time-updating scheme}''. Submitted on June 1, 2019; revised on July 23, 2019; Accepted on July 24, 2019.
}

 
\abstract{This paper deals with the concept of curvature of framed space curves, their higher-order derivatives, variations, and co-rotational derivatives. We realize that parametrizing rotation tensor using the \textit{Gibbs vector} is effective in deriving a closed-form formula to obtain any order derivative of the curvature tensor as the summation of functions of the parametrizing quantity and its derivatives. We use these results for formulating a linearized updating algorithm for curvature and its derivatives when the configuration of the curve acquires a small increment. Finally, the MATLAB code to obtain updated curvature (spatial and material) and its derivatives is presented.
}

\begin{document}
\newpage
\tableofcontents
\listoffigures
\listoftables
\newpage
\section{Introduction}
\label{section1}
In any standard course on differential geometry, smooth space curves are the first structures taught. This is generally because space curves are manifolds of a single dimension. The study of curves dates back to the $17^{\text{th}}$ century, with the primary works contributed by numerous lumenaries such as Descartes \cite{descartes1637geometrie}, Euler \cite{euler1748introductio}, Frenet \cite{frenet1852courbes}, Serret \cite{serret1851quelques}, and the contribution of Fermat mentioned in \cite{coolidge1951story}. Perhaps the most important contribution to curve framing can be attributed to Frenet-Serret \cite{frenet1852courbes}, \cite{serret1851quelques} and Bishop \cite{bishop1975there}. In the case of Frenet-Serret framing, the curve is characterized by the coordinate system-invariant scalar curvature and torsion, which require the curve to be regular, minimally $C^3$ continuous, and non-degenerate for unique existence. Bishop's frame \cite{bishop1975there}, also called the \textit{Relatively Parallel Adapted Frame}, may be used to uniquely frame a regular, minimally $C^2$ continuous curve using two invariants (known as the \textit{normal development} of the curve) that can be uniquely defined if we specify the orthogonal vectors spanning the normal plane of such a curve at a particular point on it. The two invariant quantities for Frenet-Serret frame and the normal development of Bishop's frame constitute the components of the curvature vector, also known as \textit{Darboux vector}. The curvature vector is essential in understanding the local deviation of the curve from its tangent vector.

Both of these framing techniques are intrinsic to the curve itself. The Frenet-Serret frame suffers from singularity and non-uniqueness at a point of inflection. Furthermore, these framing techniques are not the best option to describe configuration of physical systems like, for example, the configuration space of drones, the kinematics of geometrically-exact beams, the configuration of a rotating top, etc. In such cases, a general frame that does not necessarily contain the tangent vector and is defined by the configuration of the system becomes more valuable. This made us (refer to Chadha and Todd \cite{chadha2019material}) propose a \textit{material frame} useful for considering practical problems like computer graphics, trajectory tracking, path estimation, shape-sensing \cite{chadha2017generalized}, the structure of DNA, geometric design of stair cases, etc. With appropriate restrictions and constraints, any frame may be obtained from the general material frame (refer to section 4 of \cite{chadha2019material}). 

Among the numerous practical applications of framed space curves, extensive research efforts in the area of geometrically-exact beams subjected to large deformation and finite strains have been made in the past (see \cite{cosserat1909theorie}, \cite{ericksen1957exact}, \cite{simo1985finite}, \cite{chadha2017introductory} and the references therein). Seminal contributions have been made tackling theoretical and computational techniques associated with modeling and solving for geometrically exact beams, e.g., problems related to interpolation and discretization techniques (see \cite{zupan2003finite}), composite beams (see \cite{hodges2006nonlinear}) and variational principles (see \cite{simo1988hamiltonian},\cite{vetyukov2014nonlinear}). In this paper, we focus our attention on the curvature of the \textit{material frame} and its derivatives. We were presented the need for obtaining higher-order derivatives of curvature while investigating higher-order geometrically-exact beam/rod theory. The kinematics of beam/rods under arbitrarily-large deformations defined in Chadha and Todd \cite{chadha2018comprehensive} renders the deformation map not only to be a function of curvature, but also a function of its higher-order derivatives. Numerical solution of such problems using Finite Element Analysis needs updating of these kinematic quantities. In such problems, derivatives of the curvature tensor gain importance. Apart from a practical viewpoint, the fact that the Lie proper orthogonal rotational group $SO(3)$ and its Lie Algebra $so(3)$ constitute a central role in the area of Lie group theory makes it worthwhile to investigate the higher-order partial and co-rotational derivatives of curvature and the associated quantities. 

After this introductory section, this work consist of four additional sections. In section \ref{section2}, we briefly summarize concepts related to framed space curves, finite rotation parameterized by the rotation tensor, material and spatial curvature, and variation and linearization of the rotation tensor. In section \ref{section3}, we start by presenting the formula for the derivative of Lie-bracket. We present a particular way to parameterize the rotation tensor in a manner that is beneficial to obtain a closed-form formula for the derivatives of \textit{spatial} curvature. Once the derivatives of spatial curvature are obtained, we present an approach to derive the \textit{co-rotational} derivatives of curvature and finally obtain an expression for the derivative of \textit{material} curvature. Section \ref{section4} deals with updating the curvature and its derivatives given the \textit{current incremental rotation vector} and its derivatives. Section \ref{conclusion} concludes the paper.   

We conclude this introductory section with a note on notation and definitions. The $n$ dimensional Euclidean space is represented by $\mathbb{R}^n$, with $\mathbb{R}^1\equiv\mathbb{R}$. The space of real numbers and integers is denoted by $\mathbb{R}$ and $\mathbb{Z}$, with $\mathbb{R}^+$ and $\mathbb{Z}^+$ giving the set of positive real numbers and integers (including 0) respectively. The dot product, ordinary vector product and tensor product of two Euclidean vectors $\boldsymbol{v}_1$ and $\boldsymbol{v}_2$ are defined as $\boldsymbol{v}_1\cdot\boldsymbol{v}_2$, $\boldsymbol{v}_1\times\boldsymbol{v}_2$, and $\boldsymbol{v}_1\otimes\boldsymbol{v}_2$ respectively. The Euclidean norm is represented by $\left\lVert.\right\rVert$ or the un-bolded version of the symbol (for example, $\left\lVert\boldsymbol{v}\right\rVert\equiv v$). Secondly, $n^{\text{th}}$ (with $n\in\mathbb{Z}^+$) order partial derivative with respect to a scalar quantity, $\xi$ for instance, is given by the operator $\frac{\partial^n}{\partial\xi^n}=\partial^n_{\xi}$. For $n=1$, we define $\partial^1_{\xi}\equiv\partial_{\xi}$ and note that for $n=0$, $\partial^0_{\xi}$ is an identity operator. A vector, tensor or a matrix is represented by bold symbol and their components are given by indexed un-bolded symbols. For $i,j\in\mathbb{Z}^+$, the Kronecker delta function is defined as $\delta_{ij}=\begin{cases}0 & \text{ if }i\neq j\\ 1& \text{ if }i=j\end{cases}$. The action of a tensor $\boldsymbol{A}$ onto the vector $\boldsymbol{v}$ is represented by $\boldsymbol{A}\boldsymbol{v}\equiv\boldsymbol{A}.\boldsymbol{v}$. We note that the centered dot ``$\cdot$'' is meant for dot product between two vectors, whereas the action of a tensor onto the vector, the matrix multiplication or product of a scalar to a matrix (or a vector) is denoted by a lower dot ``$.$''. For $n,i\in\mathbb{Z}^+$ and $n\geq i\geq0$, the binomial coefficient is defined as $C^n_i=\frac{n!}{i!(n-i)!}$. We note two useful properties of binomial coefficient in \textit{Theorem 0}.

\paragraph{Theorem 0:} For $i,n\in\mathbb{Z}^+$ and $i\leq n$, the following holds
\begin{subequations}
	\begin{align}
	&	C_i^n=C^n_{(n-i)};\label{binomial_prop3}\\
		&C_i^n=\begin{cases}1 & \text{ if } i=0 \text{ or }i=n\\
	C_{(i-1)}^{(n-1)}+C^{(n-1)}_i & \text{ otherwise}.
	\end{cases}\label{binomial_prop1}.
	\end{align}
\end{subequations}
\paragraph{Proof:} Result \eqref{binomial_prop3} follows from the definition of binomial coefficient. The recurrence-formula \eqref{binomial_prop1} is obtained from the result $C_{(i+1)}^{(n+1)}=C_i^n+C_{(i+1)}^n$, that is easily provable using the definition of binomial coefficient.
$\square$

\section{Framed space curves}
\label{section2}
In general terms, the application of this paper is to study a smooth curve on $SO(3)$ manifold, a proper orthogonal Lie group. Physically, it is easy to understand the problem at hand by considering a framed space curve. This is because the orthogonal frame field attached to the space curve represent a curve on $SO(3)$. 
\subsection{Finite rotation}
A framed space curve, parameterized by the arc-length $\xi\in[0,L]$, is defined by the position vector $\boldsymbol{\varphi}(\xi)\in\mathbb{R}^3$ and the orthonormal material frame field $\{\boldsymbol{d}_i(\xi)\}$. Let $\{\boldsymbol{E}_i\}$ define a fixed orthonormal reference frame such that we may define the orthogonal rotation tensor $\boldsymbol{Q}(\xi)$ as:
\begin{equation} \label{Q}
\begin{aligned}
\boldsymbol{d}_i(\xi)&=\boldsymbol{Q}(\xi).\boldsymbol{E}_i;\\
\boldsymbol{Q}(\xi)&=\sum_{i=1}^3\boldsymbol{d}_i(\xi)\otimes\boldsymbol{E}_i;\\
	[\boldsymbol{Q}]_{\boldsymbol{E}_i\otimes\boldsymbol{E}_j}&=\sum_{i,j=1}^3Q_{ij}\boldsymbol{E}_i\otimes\boldsymbol{E}_j;\quad Q_{ij}=\boldsymbol{E}_i\cdot\boldsymbol{d}_j.
\end{aligned}
\end{equation}
In the equations above, $\boldsymbol{Q}(\xi)$ denotes the orthogonal rotation tensor field and $	[\boldsymbol{Q}]_{\boldsymbol{E}_i\otimes\boldsymbol{E}_j}$ represents the dyadic form of the rotation tensor with components $Q_{ij}$ defined in the $\{\boldsymbol{E}_i\}$ frame. Finite rotations are represented by an element of a proper orthogonal rotation group $SO(3)$. The $SO(3)$ manifold is a non-linear compact Lie group that has a linear skew-symmetric matrix as its Lie algebra, $so(3)$. The Lie algebra to $SO(3)$ represents its tangent plane at the identity $\boldsymbol{I}_3\in SO(3)$. The $SO(3)$ manifold and its Lie algebra $so(3)$ are defined as
\begin{subequations} 
	\begin{gather}
	SO(3):=\{\boldsymbol{Q}:\mathbb{R}^3\longrightarrow\mathbb{R}^3|\text{ } \boldsymbol{Q}^T\boldsymbol{Q}=\boldsymbol{I}_3 \text{, and } \det\boldsymbol{Q}=1\}; \label{SO3}\\
	so(3):=\{\hat{\boldsymbol\uptheta}:\mathbb{R}^3\longrightarrow\mathbb{R}^3|\text{ }\hat{\boldsymbol\uptheta}\text{ is linear, and }\hat{\boldsymbol\uptheta}+\hat{\boldsymbol\uptheta}^T=\boldsymbol{0}_3\}. \label{so3}
	\end{gather}
\end{subequations}
In the equation above, $\boldsymbol{0}_3$ represents $3\times3$ zero matrix, whereas $\boldsymbol{I}_3=\sum_{i=1}^3\boldsymbol{E}_i\otimes\boldsymbol{E}_i=\sum_{i,j=1}^3\delta_{ij}\boldsymbol{E}_i\otimes\boldsymbol{E}_j$ is the identity tensor with respect to which the director frame field is calibrated. The anti-symmetric tensor $\hat{\boldsymbol\uptheta}\in so(3)$ is equivalent to the associated axial vector $\boldsymbol\uptheta\in \mathbb{R}^3$ in the sense that for any vector $\boldsymbol{v}\in\mathbb{R}^3$, we have $\hat{\boldsymbol\uptheta}.\boldsymbol{v}=\boldsymbol\uptheta\times\boldsymbol{v}$. Therefore, there exist an isomorphism between $\mathbb{R}^3$ and $so(3)$. The action of $\hat{\boldsymbol\uptheta}$ onto the vector $\boldsymbol{v}$ (yielding $\hat{\boldsymbol\uptheta}.\boldsymbol{v}$) results into an infinitesimal rotation of the vector $\boldsymbol{v}$ about the unit vector $\frac{\boldsymbol\uptheta}{\left\lVert\boldsymbol\uptheta\right\rVert}$ by an amount $\left\lVert\boldsymbol\uptheta\right\rVert$ (referred to as rotation about $\boldsymbol\uptheta$ vector; hence, $\boldsymbol\uptheta$ is called an \textit{axial vector}). From here on, any matrix quantity with a hat on it ($\hat.$) represents an anti-symmetric matrix. For later use, we define zero vector as $\boldsymbol{0}_1=[0,0,0]^T$. At this point, we define the Lie-bracket of two anti-symmetric matrix as $[.,.]:so(3)\times so(3)\longrightarrow \mathbb{R}^3$, such that for any $\hat{\boldsymbol{a}},\hat{\boldsymbol{b}}\in so(3)$ with corresponding axial vectors $\boldsymbol{a},\boldsymbol{b}\in \mathbb{R}^3$ respectively and any vector $\boldsymbol{v}\in\mathbb{R}^3$, we have
\begin{subequations}\label{lie_brackets}
\begin{gather}
\left[\hat{\boldsymbol{a}},\hat{\boldsymbol{b}}\right]=(\hat{\boldsymbol{a}}.\hat{\boldsymbol{b}}-\hat{\boldsymbol{b}}.\hat{\boldsymbol{a}});\label{lie_brackets1}\\\left[\hat{\boldsymbol{a}},\hat{\boldsymbol{b}}\right].\boldsymbol{v}=(\boldsymbol{a}\times\boldsymbol{b})\times\boldsymbol{v}.\label{lie_brackets2}
\end{gather}
\end{subequations}
We note two important properties of Lie-bracket:
\begin{subequations}\label{lie_brackets_properties}
	\begin{gather}
	\left[\hat{\boldsymbol{a}},\hat{\boldsymbol{a}}\right]=\boldsymbol{0}_3;\label{lie_brackets_properties1}\\\left[\hat{\boldsymbol{a}},\hat{\boldsymbol{b}}\right]=-\left[\hat{\boldsymbol{b}},\hat{\boldsymbol{a}}\right].\label{lie_brackets_properties2}
	\end{gather}
\end{subequations}

The definition of $SO(3)$ in Eq. \eqref{SO3} allows rotation tensor to be parameterized by a rotation vector $\boldsymbol{\uptheta}\in\mathbb{R}^3$ (with corresponding anti-symmetric matrix $\hat{\boldsymbol\uptheta}\in so(3)$). This is attained by means of exponential map $\text{exp}:so(3)\longrightarrow SO(3)$ such that,
\begin{subequations}
	\begin{gather}
	\boldsymbol{Q}(\boldsymbol{\uptheta})=\boldsymbol{I}_3+\frac{\sin\uptheta}{{\uptheta}} \ \hat{\boldsymbol\uptheta}+\frac{(1-\cos\uptheta)}{{\uptheta^2}}\hat{\boldsymbol\uptheta}^2   =\sum_{i\geq0}\frac{\hat{\boldsymbol\uptheta}^i}{i!}=\text{exp}(\hat{\boldsymbol\uptheta});\label{exp_map1}\\
	\boldsymbol{Q}^T(\boldsymbol{\uptheta})=\boldsymbol{Q}^{-1}(\boldsymbol{\uptheta})=\boldsymbol{I}_3-\frac{\sin\uptheta}{{\color{black}{\uptheta}}} \ \hat{\boldsymbol\uptheta}+\frac{(1-\cos\theta)}{{\uptheta^2}}\hat{\boldsymbol\uptheta}^2   =\text{exp}(-\hat{\boldsymbol\uptheta}).\label{exp_map2}
	\end{gather}
\end{subequations}
Here $\hat{\boldsymbol\uptheta}^0=\boldsymbol{I}_3$. Subtracting Eq. \eqref{exp_map2} from \eqref{exp_map1}, we obtain the associated anti-symmetric matrix $\hat{\boldsymbol\uptheta}$ as
\begin{equation}\label{uptheta_hat}
\hat{\boldsymbol{\uptheta}}=\frac{\uptheta}{2 \sin{\uptheta}}(\boldsymbol{Q}-\boldsymbol{Q}^T).
\end{equation}
The exponential map is a local homeomorphism in the neighborhood of identity $\boldsymbol{I}_3\in SO(3)$ for $\uptheta\in[0,\pi)$. Refer to section 3.2.2 of Chadha and Todd \cite{chadha2019material} for further discussion. The local homeomorphism of $\text{exp}$ map guarantees the existence of a unique inverse of exponential map in the neighborhood of $\boldsymbol{I}_3\in SO(3)$, called the logarithm map $\text{log}:SO(3)\longrightarrow so(3)$, such that
\begin{equation}\label{log}
\text{log}\left(\boldsymbol{Q}(\boldsymbol{\uptheta})\right)=\text{log}(\text{exp}(\hat{\boldsymbol\uptheta}))=\hat{\boldsymbol\uptheta}\in so(3).
\end{equation}
The norm of logarithm map is defined as the Euclidean norm of the associated rotation vector as
\begin{equation}\label{S321_Eq4d}
\begin{gathered}
\|\text{log}\left(\boldsymbol{Q}(\boldsymbol{\uptheta})\right)\|=\uptheta=\sqrt{\frac{1}{2}\text{Tr}\left(\hat{\boldsymbol\uptheta}.\hat{\boldsymbol\uptheta}\right)}.
\end{gathered}
\end{equation}
Equation \ref{S321_Eq4d} above defines a metric that is useful in measuring errors in the director triad. The rotation tensor can also be represented by means of unit quaternion. For a detailed discussion on finite rotation, refer to Argyris \cite{argyris1982excursion}, Ibrahimbegovic \cite{ibrahimbegovic1995computational} and Diebel \cite{diebel2006representing}.
\subsection{Curvature: Material and spatial form}
Curvature tensor defines local change of the triad. We have
\begin{equation}\label{curvature}
\partial_{\xi}\boldsymbol{d}_i=\partial_{\xi}\boldsymbol{Q}.\boldsymbol{E}_i=\partial_{\xi}\boldsymbol{Q}.\boldsymbol{Q}^T.\boldsymbol{d}_i=\hat{\boldsymbol{\kappa}}.\boldsymbol{d}_i.
\end{equation}
Here, $\hat{\boldsymbol{\kappa}}=\partial_{\xi}\boldsymbol{Q}.\boldsymbol{Q}^T$ gives curvature tensor. It is an anti-symmetric matrix with the corresponding axial vector $\boldsymbol{\kappa}=\overline\kappa_i\boldsymbol{d}_i$, known as curvature vector. We define $T_{\boldsymbol{Q}}SO(3)$ as the tangent plane of non-linear $SO(3)$ manifold, such that $\partial_{\xi}\boldsymbol{Q}=\hat{\boldsymbol{\kappa}}.\boldsymbol{Q}\in T_{\boldsymbol{Q}}SO(3)$. We note that $so(3)= T_{\boldsymbol{I}_3}SO(3)$. Here, we can interpret the director field $\boldsymbol{Q}(\xi)$ as a curve on $SO(3)$ parameterized by $\xi\in[0,L]$. Therefore, $\partial_{\xi}\boldsymbol{Q}\in T_{\boldsymbol{Q}}SO(3)$ gives the tangent vector field to the curve $\boldsymbol{Q}(\xi)$. If rotation tensor $\boldsymbol{Q}$ is parameterized by rotation vector $\boldsymbol{\uptheta}$, using Eq. \eqref{Q} and \eqref{curvature}, we obtain the following expression of the curvature tensor $\hat{\boldsymbol{\kappa}}$ as
\begin{equation}\label{curv_tensor_uptheta}
\hat{\boldsymbol{\kappa}}=\partial_{\xi}\text{exp}(\hat{\boldsymbol{\uptheta}}).\text{exp}(-\hat{\boldsymbol{\uptheta}})=\left(\frac{\sin{\uptheta}}{\uptheta}\right)\partial_{\xi}\hat{\boldsymbol{\uptheta}}+\left(\frac{1-\cos{\uptheta}}{\uptheta^2}\right)\left[\hat{\boldsymbol{\uptheta}},\partial_{\xi}\hat{\boldsymbol{\uptheta}}\right]+(\boldsymbol\uptheta\cdot\partial_{\xi}\boldsymbol\uptheta)\left(\frac{\uptheta-\sin{\uptheta}}{\uptheta^3}\right)\hat{\boldsymbol{\uptheta}}.
\end{equation}
Using equation set \eqref{lie_brackets} and the result in Eq. \eqref{curv_tensor_uptheta}, we obtain the curvature vector $\boldsymbol{\kappa}$ as
\begin{equation}\label{curv_vec_uptheta}
\begin{gathered}
\boldsymbol{\kappa}=\boldsymbol{T}_{\uptheta}.\partial_{\xi}\boldsymbol{\uptheta};\\
\boldsymbol{T}_{\uptheta}=\left(\frac{\sin{\uptheta}}{\uptheta}\right)\boldsymbol{I}_3+\left(\frac{1-\cos{\uptheta}}{\uptheta^2}\right)\hat{\boldsymbol\uptheta}+\left(\frac{\uptheta-\sin{\uptheta}}{\uptheta^3}\right)\boldsymbol{\uptheta}\otimes\boldsymbol{\uptheta}.
\end{gathered}
\end{equation}
Similarly, we have
\begin{equation}\label{uptheta_curv_vec}
\begin{gathered}
\partial_{\xi}\boldsymbol{\uptheta}=\boldsymbol{T}_{\uptheta}^{-1}.\boldsymbol{\kappa};\\
\boldsymbol{T}_{\uptheta}^{-1}=\left(\frac{1}{2}\frac{\uptheta}{\tan{\frac{\uptheta}{2}}}\right)\boldsymbol{I}_3-\frac{\hat{\boldsymbol\uptheta}}{2}+\frac{1}{\uptheta^2}\left(1-\frac{1}{2}\frac{\uptheta}{\tan{\frac{\uptheta}{2}}}\right)\boldsymbol{\uptheta}\otimes\boldsymbol{\uptheta}.
\end{gathered}
\end{equation}
Readers may refer to Ibrahimbegovic et al. \cite{ibrahimbegovic1995computational} for the derivation of $\boldsymbol{T}_{\uptheta}^{-1}$. 

It is interesting to interpret the curvature vector $\boldsymbol{\kappa}$ and the derivative of the rotation vector $\partial_{\xi}\boldsymbol\uptheta$ from a physical viewpoint. At an arc-length $\xi$, the director triad $\{\boldsymbol{d}_i(\xi)\}$ rotates about the vector $\boldsymbol{\kappa}(\xi).\text{d}\xi$ to yield the triad at $\{\boldsymbol{d}_i(\xi+\text{d}\xi)\}$. The triad $\{\boldsymbol{d}_i(\xi)\}$ and $\{\boldsymbol{d}_i(\xi+\text{d}\xi)\}$ are obtained by finite rotation of the frame $\{\boldsymbol{E}_i\}$ about the rotation vector $\boldsymbol{\uptheta}(\xi)$ and $\boldsymbol{\uptheta}(\xi+\text{d}\xi)=\boldsymbol{\uptheta}(\xi)+\partial_{\xi}\uptheta(\xi).\text{d}\xi$ respectively. In terms of the exponential map,
\begin{equation}\label{Q4}
\begin{gathered}
\boldsymbol{Q}(\xi+\text{d}\xi)=\text{exp}(\hat{\boldsymbol{\kappa}}(\xi).\text{d}\xi).\boldsymbol{Q}(\xi)=\text{exp}(\hat{\boldsymbol{\kappa}}(\xi).\text{d}\xi).\text{exp}(\hat{\boldsymbol{\uptheta}}(\xi));\\
\boldsymbol{Q}(\xi+\text{d}\xi)=\boldsymbol{Q}(\boldsymbol\uptheta(\xi+\text{d}\xi))=\text{exp}(\hat{\boldsymbol\uptheta}(\xi)+\partial_{\xi}\hat{\boldsymbol{\uptheta}}(\xi).\text{d}\xi).
\end{gathered}
\end{equation}

We define the material curvature tensor as $\hat{\overline{\boldsymbol{\kappa}}}=\boldsymbol{Q}^T.\hat{\boldsymbol{\kappa}}.\boldsymbol{Q}=\boldsymbol{Q}^T.\partial_{\xi}\boldsymbol{Q}\in so(3)$ obtained by parallel transport of $\hat{\boldsymbol{\kappa}}.\boldsymbol{Q}$ from $T_{\boldsymbol{Q}}SO(3)\longrightarrow so(3)$. From Eq. \eqref{Q}, $\boldsymbol{Q}$ represents the finite rotation, while $\hat{\overline{\boldsymbol{\kappa}}}$ gives an infinitesimal rotation with respect to the calibrating frame (fixed reference frame) $\{\boldsymbol{E}_i\}$, and $\boldsymbol{Q}.\hat{\overline{\boldsymbol{\kappa}}}=\hat{\boldsymbol{\kappa}}.\boldsymbol{Q}$ represents an infinitesimal rotation with respect to the $\{\boldsymbol{d}_i\}$ frame. In the physical context of rotation, the tangent vector $\boldsymbol{Q}.\hat{\overline{\boldsymbol{\kappa}}}$ and $\boldsymbol{\kappa}.\boldsymbol{Q}$ performs an infinitesimal rotation with respect to $\{\boldsymbol{d}_i\}$ frame, but the quantity $\boldsymbol{Q}.\hat{\overline{\boldsymbol{\kappa}}}$ is obtained by a left translation of $\hat{\overline{\boldsymbol{\kappa}}}\in so(3)$ to $\boldsymbol{Q}.\hat{\overline{\boldsymbol{\kappa}}}\in T_{\boldsymbol{Q}}SO(3)$. Moreover, $\hat{\boldsymbol{\kappa}}.\boldsymbol{Q}$ represents the superimposition of infinitesimal rotation contributed by $\hat{\boldsymbol{\kappa}}$ onto the finite rotation due to $\boldsymbol{Q}$ (this is also called the right translation of $\hat{\boldsymbol{\kappa}}\in so(3)$ to the tangent vector $\hat{\boldsymbol{\kappa}}.\boldsymbol{Q}\in T_{\boldsymbol{Q}}SO(3)$). The former kind of tangent vector fields are known as \textit{left-invariant} and the latter as\textit{ right-invariant} fields. We also observe that
\begin{equation}
\begin{gathered}
\left[\hat{\boldsymbol{\kappa}}\right]_{\boldsymbol{d}_i\otimes\boldsymbol{d}_j}=\left[\hat{\boldsymbol{\overline \kappa}}\right]_{\boldsymbol{E}_i\otimes\boldsymbol{E}_j}=\begin{bmatrix}
0 & -\overline{\kappa}_3 &  \overline{\kappa}_2 \\
\overline{\kappa}_3 & 0 &  -\overline{\kappa}_1 \\
-\overline{\kappa}_2 & \overline{\kappa}_1 &  0
\end{bmatrix}.
\end{gathered}
\end{equation}
Let $\boldsymbol{\kappa}$ and $\overline{\boldsymbol{\kappa}}$ represent the axial vector corresponding to the anti-symmetric matrix $\hat{\boldsymbol{\kappa}}$ and $\hat{\overline{\boldsymbol{\kappa}}}$ respectively. It can then be proved that $\overline{\boldsymbol{\kappa}}=\boldsymbol{Q}^T.\boldsymbol{\kappa}$ such that if $\boldsymbol{\kappa}=\overline{\kappa}_i\boldsymbol{d}_i$, then $\overline{\boldsymbol{\kappa}}=\overline{\kappa}_i\boldsymbol{E}_i$.  We call the quantities $\hat{\overline{\boldsymbol{\kappa}}}$ and $\overline{\boldsymbol{\kappa}}$ as \textit{material representation}; and  $\hat{\boldsymbol{\kappa}}$ and $\boldsymbol{\kappa}$ as \textit{spatial representation} of the \textit{curvature tensor} and the \textit{curvature vector}, respectively. 
\subsection{Material and spatial quantities and their derivatives}
As we have defined the material and spatial forms of the curvature vector (and tensor), it is rather useful to define a vector $\boldsymbol{v}\in\mathbb{R}^3$ in its material and spatial form. Consider a spatial and material vector $\boldsymbol{v}=\overline{v}_i\boldsymbol{d}_i$ and $\overline{\boldsymbol{v}}=\overline{v}_i\boldsymbol{E}_i$, respectively, such that $\boldsymbol{v}=\boldsymbol{Q}.\overline{\boldsymbol{v}}$. The derivative of these vectors are obtained as
\begin{equation}\label{corotational1}
\begin{gathered}
\partial_{\xi}\boldsymbol{v}=\partial_{\xi}\overline{v}_i.\boldsymbol{d}_i+\overline{v}_i.\partial_{\xi}\boldsymbol{d}_i=\tilde\partial_{\xi}\boldsymbol{v}+\boldsymbol{\kappa}\times\boldsymbol{v};\\
\partial_{\xi}\overline{\boldsymbol{v}}=\partial_{\xi}\overline{v}_i.\boldsymbol{E}_i=\boldsymbol{Q}^T.\tilde\partial_{\xi}\boldsymbol{v}.
\end{gathered}
\end{equation}
In the equation above, $\tilde\partial_{\xi}\boldsymbol{v}$ defines co-rotational derivative of spatial vector $\boldsymbol{v}$. It essentially gives the change in components of the vector $\boldsymbol{v}$, provided the frame of reference is assumed to be fixed. Geometrically, it is obtained by parallel-transport (left translation) of the vector $\partial_{\xi}\overline{\boldsymbol{v}}$. 

Along similar lines, consider a spatial and material tensor $\boldsymbol{A}=\sum_{i,j=1}^3\overline{A}_{ij}\boldsymbol{d}_i\otimes\boldsymbol{d}_j$ and $\overline{\boldsymbol{A}}=\sum_{i,j=1}^3\overline{A}_{ij}\boldsymbol{E}_i\otimes\boldsymbol{E}_j$ respectively, such that $\boldsymbol{A}=\boldsymbol{Q}.\overline{\boldsymbol{A}}.\boldsymbol{Q}^T$ (or equivalently, $\overline{\boldsymbol{A}}=\boldsymbol{Q}^T.\boldsymbol{A}.\boldsymbol{Q}$). Realizing $\partial_{\xi}\boldsymbol{Q}=\hat{\boldsymbol{\kappa}}.\boldsymbol{Q}$ and $\partial_{\xi}\boldsymbol{Q}^T=-\boldsymbol{Q}^T.\hat{\boldsymbol{\kappa}}$, we have the following
\begin{equation}\label{corotational2}
\begin{aligned}
\partial_{\xi}\boldsymbol{A}&=\boldsymbol{Q}.\partial_{\xi}\overline{\boldsymbol{A}}.\boldsymbol{Q}^T+\partial_{\xi}\boldsymbol{Q}.\overline{\boldsymbol{A}}.\boldsymbol{Q}^T+\boldsymbol{Q}.\overline{\boldsymbol{A}}.\partial_{\xi}\boldsymbol{Q}^T\\
&=\boldsymbol{Q}.\partial_{\xi}\overline{\boldsymbol{A}}.\boldsymbol{Q}^T+\hat{\boldsymbol{\kappa}}.(\boldsymbol{Q}.\overline{\boldsymbol{A}}.\boldsymbol{Q}^T)-(\boldsymbol{Q}.\overline{\boldsymbol{A}}.\boldsymbol{Q}^T).\hat{\boldsymbol{\kappa}}\\
&=\tilde\partial_{\xi}\boldsymbol{A}+\hat{\boldsymbol{\kappa}}.\boldsymbol{A}-\boldsymbol{A}.\hat{\boldsymbol{\kappa}}.
\end{aligned}
\end{equation}
In the equation above, we define the co-rotational derivative of the tensor $\boldsymbol{A}$ as $\tilde\partial_{\xi}\boldsymbol{A}=\boldsymbol{Q}.\partial_{\xi}\overline{\boldsymbol{A}}.\boldsymbol{Q}^T$. Physically, it gives the change in components of tensor $\boldsymbol{A}$ setting the reference frame constant. From here on, $\tilde\partial_{x}(.)$ represents the co-rotational derivative of some quantity $(.)$ with respect to the variable $x$. We now present some propositions describing higher order co-rotational derivatives. 
\paragraph{Proposition 1:} For any vector $\boldsymbol{v}\in\mathbb R^3$, define the operator $\hat\partial_{\xi}$ such that $\hat\partial_{\xi}^n\boldsymbol{v}=\hat{\boldsymbol{\kappa}}^n.\boldsymbol{v}$, where, for example $\hat{\boldsymbol{\kappa}}^3=\hat{\boldsymbol{\kappa}}.\hat{\boldsymbol{\kappa}}.\hat{\boldsymbol{\kappa}}$. The $n^\text{th}$ order co-rotational derivative $\tilde\partial_{\xi}^n$ is then given by $(\partial_{\xi}-\hat{\partial}_{\xi})^n$ such that the order of operations are not commutative (for example, $\partial_{\xi}\hat{\partial}_{\xi}\neq\hat{\partial}_{\xi}\partial_{\xi}$) and $\tilde\partial_{\xi}^0=(\partial_{\xi}-\hat{\partial}_{\xi})^0$ is an identity operator.

\paragraph{Proof:} We prove proposition 1 by using the principle of mathematical induction. Consider a vector $\boldsymbol{v}\in\mathbb R^3$, assuming that $\tilde\partial_{\xi}^n\boldsymbol{v}=(\partial_{\xi}-\hat{\partial}_{\xi})^n\boldsymbol{v}$, for $n=1,2$ and $3$, we have 
\begin{equation}\label{corot_der1}
\begin{aligned}
&\tilde\partial_{\xi}^1\equiv\tilde\partial_{\xi}=(\partial_{\xi}-\hat{\partial}_{\xi});\\
&\tilde\partial_{\xi}^2=(\partial_{\xi}-\hat{\partial}_{\xi})^2=\partial_{\xi}^2-\partial_{\xi}\hat{\partial}_{\xi}-\hat{\partial}_{\xi}\partial_{\xi}+\hat{\partial}_{\xi}^2;\\
&\tilde\partial_{\xi}^3=(\partial_{\xi}-\hat{\partial}_{\xi})^3=\partial_{\xi}^3-\partial_{\xi}^2\hat{\partial}_{\xi}-\partial_{\xi}\hat{\partial}_{\xi}\partial_{\xi}+\partial_{\xi}\hat{\partial}_{\xi}^2-\hat{\partial}_{\xi}\partial_{\xi}^2+\hat{\partial}_{\xi}^2\partial_{\xi}+\hat{\partial}_{\xi}\partial_{\xi}\hat{\partial}_{\xi}-\hat{\partial}_{\xi}^3;
\end{aligned}
\end{equation}
such that
\begin{subequations}\label{corot_der}
	\begin{align}
	&\tilde\partial_{\xi}\boldsymbol{v}=\partial_{\xi}\boldsymbol{v}-\hat{\partial}_{\xi}\boldsymbol{v}=\partial_{\xi}\boldsymbol{v}-\hat{\boldsymbol\kappa}.\boldsymbol{v};\label{corot_der2}\\
	&\begin{aligned}
	\tilde\partial_{\xi}^2\boldsymbol{v}&=\partial_{\xi}^2\boldsymbol{v}-\partial_{\xi}\hat{\partial}_{\xi}\boldsymbol{v}-\hat{\partial}_{\xi}\partial_{\xi}\boldsymbol{v}+\hat{\partial}_{\xi}^2\boldsymbol{v}=\partial_{\xi}^2\boldsymbol{v}-\partial_{\xi}(\hat{\boldsymbol\kappa}.\boldsymbol{v})-\hat{\boldsymbol\kappa}.\partial_{\xi}\boldsymbol{v}+\hat{\boldsymbol\kappa}.\hat{\boldsymbol\kappa}.\boldsymbol{v}\\&=\partial_{\xi}^2\boldsymbol{v}+(\hat{\boldsymbol\kappa}.\hat{\boldsymbol\kappa}-\partial_{\xi}\hat{\boldsymbol\kappa}).\boldsymbol{v}-2\hat{\boldsymbol\kappa}.\partial_{\xi}\boldsymbol{v};
	\end{aligned}\label{corot_der3}\\
	&\begin{aligned}
	\tilde\partial_{\xi}^3\boldsymbol{v}&=\partial_{\xi}^3\boldsymbol{v}-\partial_{\xi}^2(\hat{\boldsymbol\kappa}.\boldsymbol{v})-\partial_{\xi}(\hat{\boldsymbol\kappa}.\partial_{\xi}\boldsymbol v)+\partial_{\xi}(\hat{\boldsymbol\kappa}.\hat{\boldsymbol\kappa}.\boldsymbol{v})-\hat{\boldsymbol\kappa}.\partial_{\xi}^2\boldsymbol{v}+\hat{\boldsymbol\kappa}.\hat{\boldsymbol\kappa}.\partial_{\xi}\boldsymbol{v}+\hat{\partial}_{\xi}\partial_{\xi}(\hat{\boldsymbol\kappa}.\boldsymbol{v})\\&\quad-\hat{\boldsymbol\kappa}.\hat{\boldsymbol\kappa}.\hat{\boldsymbol\kappa}.\boldsymbol{v}=\partial_{\xi}^3\boldsymbol{v}-3\hat{\boldsymbol\kappa}.\partial_{\xi}^2\boldsymbol{v}+(\partial_{\xi}\hat{\boldsymbol\kappa}.\hat{\boldsymbol\kappa}+2\hat{\boldsymbol\kappa}.\partial_{\xi}\hat{\boldsymbol\kappa}-\partial_{\xi}^2\hat{\boldsymbol\kappa}-\hat{\boldsymbol\kappa}.\hat{\boldsymbol\kappa}.\hat{\boldsymbol\kappa}).\boldsymbol{v}\\&\quad+(3\hat{\boldsymbol\kappa}.\hat{\boldsymbol\kappa}-3\partial_{\xi}\hat{\boldsymbol\kappa}).\partial_{\xi}\boldsymbol{v}.
	\end{aligned}\label{corot_der4}
	\end{align}
\end{subequations}
We now prove that the equation set \eqref{corot_der} may be derived using the definition of co-rotational derivatives in Eq. \eqref{corotational1}. Equation \eqref{corot_der2} is true by definition (refer to Eq. \eqref{corotational1}). Taking the derivative of Eq. \eqref{corot_der2} yields
\begin{equation}\label{corot_der5}
\begin{aligned}
&\partial_{\xi}\tilde\partial_{\xi}\boldsymbol{v}=\partial_{\xi}^2\boldsymbol{v}-\partial_{\xi}\hat{\boldsymbol\kappa}.\boldsymbol{v}-\hat{\boldsymbol\kappa}.\partial_{\xi}\boldsymbol{v}
\end{aligned}
\end{equation}
We note from the definition of co-rotational derivative in Eq. \eqref{corotational1} that
\begin{equation}\label{corot_der6}
\begin{aligned}
\partial_{\xi}\tilde\partial_{\xi}\boldsymbol{v}&=\tilde\partial_{\xi}^2\boldsymbol{v}+\hat{\boldsymbol\kappa}.\tilde\partial_{\xi}\boldsymbol{v}=\tilde\partial_{\xi}^2\boldsymbol{v}+\hat{\boldsymbol\kappa}.(\partial_{\xi}\boldsymbol{v}-\hat{\boldsymbol\kappa}.\boldsymbol{v})=\tilde\partial_{\xi}^2\boldsymbol{v}+\hat{\boldsymbol\kappa}.\partial_{\xi}\boldsymbol{v}-\hat{\boldsymbol\kappa}.\hat{\boldsymbol\kappa}.\boldsymbol{v}.
\end{aligned}
\end{equation}
From Eq. \eqref{corot_der5} and \eqref{corot_der6}, we have
\begin{equation}\label{corot_der7}
\begin{aligned}
\tilde\partial_{\xi}^2\boldsymbol{v}&=\partial_{\xi}\tilde\partial_{\xi}\boldsymbol{v}-\hat{\boldsymbol\kappa}.\partial_{\xi}\boldsymbol{v}+\hat{\boldsymbol\kappa}.\hat{\boldsymbol\kappa}.\boldsymbol{v}=\partial_{\xi}^2\boldsymbol{v}+(\hat{\boldsymbol\kappa}.\hat{\boldsymbol\kappa}-\partial_{\xi}\hat{\boldsymbol\kappa}).\boldsymbol{v}-2\hat{\boldsymbol\kappa}.\partial_{\xi}\boldsymbol{v}.
\end{aligned}
\end{equation}
Expression obtained in \eqref{corot_der7} is same as \eqref{corot_der3}. Similarly, to derive the expression for $\tilde\partial_{\xi}^3\boldsymbol{v}$, we consider
\begin{equation}\label{corot_der8}
\begin{aligned}
\partial_{\xi}\tilde\partial_{\xi}^2\boldsymbol{v}&=\partial_{\xi}^3\boldsymbol{v}-2\hat{\boldsymbol\kappa}.\partial_{\xi}^2\boldsymbol{v}+(\hat{\boldsymbol\kappa}.\hat{\boldsymbol\kappa}-3\partial_{\xi}\hat{\boldsymbol\kappa}).\partial_{\xi}\boldsymbol{v}+(\partial_{\xi}\hat{\boldsymbol\kappa}.\hat{\boldsymbol\kappa}+\hat{\boldsymbol\kappa}.\partial_{\xi}\hat{\boldsymbol\kappa}-\partial_{\xi}^2\hat{\boldsymbol\kappa}).\boldsymbol{v}.
\end{aligned}
\end{equation}
From \eqref{corotational1}, we have
\begin{equation}\label{corot_der9}
\begin{aligned}
\tilde\partial_{\xi}^3\boldsymbol{v}=\partial_{\xi}\tilde\partial_{\xi}^2\boldsymbol{v}-\hat{\boldsymbol{\kappa}}.(\partial_{\xi}\tilde\partial_{\xi}^2\boldsymbol{v}).
\end{aligned}
\end{equation}
Using the results in Eq. \eqref{corot_der8} and \eqref{corot_der9}, we arrive at the expression of $\tilde\partial_{\xi}^3\boldsymbol{v}$ as obtained in \eqref{corot_der4}. We can continue the process explained above and realize that for any $n$, $\tilde\partial_{\xi}^n=(\partial_{\xi}-\hat{\partial}_{\xi})(\partial_{\xi}-\hat{\partial}_{\xi})^{(n-1)}=(\partial_{\xi}-\hat{\partial}_{\xi})^n$. Using binomial theorem, we can also write
\begin{equation}\label{corot_der11}
\begin{aligned}
(\partial_{\xi}-\hat{\partial}_{\xi})^n=\sum_{i=0}^{n}(-1)^{(n-i)}C_i^n\partial_{\xi}^n\hat{\partial}_{\xi}^{(n-i)}.
\end{aligned}
\end{equation}
This completes the proof. $\square$
\paragraph{Proposition 2:} For any $\hat{\boldsymbol{A}}\in so(3)$ with the corresponding axial vector $\boldsymbol{A}\in\mathbb{R}^3$, the recurrence formula for the $n^{th}$ order co-rotational derivative $\tilde\partial_{\xi}^n\hat{\boldsymbol{A}}\in so(3)$ and $\tilde\partial_{\xi}^n\boldsymbol{A}\in\mathbb{R}^3$ is given as
\begin{subequations}\label{p2_1}
\begin{align}
\tilde\partial_{\xi}^n\hat{\boldsymbol{A}}&=\partial_{\xi}^n\hat{\boldsymbol{A}}-(1-\delta_{n0})\sum_{i=1}^n\partial_{\xi}^{(i-1)}\left[\hat{\boldsymbol{\kappa}},\tilde\partial_{\xi}^{(n-i)}\hat{\boldsymbol A}\right];\label{p2_1a}\\
\tilde\partial_{\xi}^n{\boldsymbol{A}}&=\partial_{\xi}^n{\boldsymbol{A}}-(1-\delta_{n0})\sum_{i=1}^n\partial_{\xi}^{(i-1)}({\boldsymbol{\kappa}}\times\tilde\partial_{\xi}^{(n-i)}{\boldsymbol A}).\label{p2_1b}
\end{align}
\end{subequations}
\paragraph{Proof:} From definition of co-rotational derivatives and Lie-bracket in Eq. \eqref{corotational2} and \eqref{lie_brackets1} respectively, we have
\begin{equation}\label{p2_2}
\begin{aligned}
\tilde\partial_{\xi}\hat{\boldsymbol A}=\partial_{\xi}\hat{\boldsymbol A}-\hat{\boldsymbol\kappa}.\hat{\boldsymbol A}+\hat{\boldsymbol A}.\hat{\boldsymbol\kappa}=\partial_{\xi}\hat{\boldsymbol A}-\left[\hat{\boldsymbol\kappa},\hat{\boldsymbol A}\right]
\end{aligned}
\end{equation}
Since, $\partial_{\xi}^m\hat{\boldsymbol A}\in so(3)$ for any $m\in\mathbb{Z}^+$, the result above can be used to obtain the recurrence-relation for $n^{\text{th}}$ order co-rotational derivative. For $n\in\mathbb{Z}^+-\{0\}$, we have
\begin{equation}\label{p2_3}
\begin{aligned}
\tilde\partial_{\xi}^n\hat{\boldsymbol{A}}=&\partial_{\xi}\left(\tilde\partial_{\xi}^{(n-1)}\hat{\boldsymbol{A}}\right)-\left[\hat{\boldsymbol{\kappa}},\tilde\partial_{\xi}^{(n-1)}\hat{\boldsymbol{A}}\right]\\
=&\partial_{\xi}\left(\partial_{\xi}\left(\tilde\partial_{\xi}^{(n-2)}\hat{\boldsymbol{A}}\right)-\left[\hat{\boldsymbol{\kappa}},\tilde\partial_{\xi}^{(n-2)}\hat{\boldsymbol{A}}\right]\right)-\left[\hat{\boldsymbol{\kappa}},\tilde\partial_{\xi}^{(n-1)}\hat{\boldsymbol{A}}\right]\\
=&\partial_{\xi}^2\left(\partial_{\xi}\left(\tilde\partial_{\xi}^{(n-3)}\hat{\boldsymbol{A}}\right)-\left[\hat{\boldsymbol{\kappa}},\tilde\partial_{\xi}^{(n-3)}\hat{\boldsymbol{A}}\right]\right)-\partial_{\xi}\left[\hat{\boldsymbol{\kappa}},\tilde\partial_{\xi}^{(n-2)}\hat{\boldsymbol{A}}\right]-\left[\hat{\boldsymbol{\kappa}},\tilde\partial_{\xi}^{(n-1)}\hat{\boldsymbol{A}}\right]\\
=&\partial_{\xi}^n\hat{\boldsymbol{A}}-\left[\hat{\boldsymbol{\kappa}},\tilde\partial_{\xi}^{(n-1)}\hat{\boldsymbol{A}}\right]-\partial_{\xi}\left[\hat{\boldsymbol{\kappa}},\tilde\partial_{\xi}^{(n-2)}\hat{\boldsymbol{A}}\right]-...-\partial_{\xi}^{(n-1)}\left[\hat{\boldsymbol{\kappa}},\hat{\boldsymbol{A}}\right]\\
=&\partial_{\xi}^n\hat{\boldsymbol{A}}-\sum_{i=1}^n\partial_{\xi}^{(i-1)}\left[\hat{\boldsymbol{\kappa}},\tilde\partial_{\xi}^{(n-i)}\hat{\boldsymbol A}\right].
\end{aligned}
\end{equation}
We note that for $n=0$, we have $\tilde\partial_{\xi}^0\hat{\boldsymbol A}=\partial_{\xi}^0\hat{\boldsymbol A}=\hat{\boldsymbol A}$. Thus, the sum part in the equation above vanish for $n=0$, justifying the use of $(1-\delta_{n0})$ factor in Eq. \eqref{p2_1a}. Result \eqref{p2_1b} follows from above derivation using Eq. \eqref{lie_brackets2}. This completes the proof. $\square$
\paragraph{Corollary 1:} The \textit{Proposition 2} can be extended for any tensor $\boldsymbol{B}$ (not necessarily an element of $so(3)$) as:
\begin{equation}\label{c1_1}
\begin{aligned}
\tilde\partial_{\xi}^n{\boldsymbol{B}}=\partial_{\xi}^n{\boldsymbol{B}}-(1-\delta_{n0})\sum_{i=1}^n\partial_{\xi}^{(i-1)}\left(\hat{\boldsymbol{\kappa}}.\tilde\partial_{\xi}^{(n-i)}{\boldsymbol B}-\tilde\partial_{\xi}^{(n-i)}{\boldsymbol B}.\hat{\boldsymbol{\kappa}}\right)
\end{aligned}
\end{equation}
\paragraph{Proof:}  This extension is possible because Lie-brackets follow chain-rule just like product of two scalar or dot product except for the fact that Lie-brackets are non-commutative (which is a stronger condition) as observed in Eq. \eqref{lie_brackets3}. $\square$
\paragraph{Proposition 3:} For spatial vector and tensor $\boldsymbol{v}$ and $\boldsymbol{A}$ respectively, with corresponding material quantities $\overline{\boldsymbol{v}}$ and $\overline{\boldsymbol{A}}$, the $n^{\text{th}}$ order co-rotational derivative $\tilde\partial_{\xi}^n\boldsymbol{v}$ and $\tilde\partial_{\xi}^n\boldsymbol{A}$  can be obtained by left-translation of the $n^{\text{th}}$ order derivative of the respective material quantities such that
\begin{subequations}\label{corot_der10}
	\begin{align}
&\tilde\partial_{\xi}^n\boldsymbol{v}=\boldsymbol{Q}.\partial_{\xi}^n\overline{\boldsymbol{v}}\\
&\tilde\partial_{\xi}^n\boldsymbol{A}=\boldsymbol{Q}.\partial_{\xi}^n\overline{\boldsymbol{A}}.\boldsymbol{Q}^T.
\end{align}
\end{subequations}
\paragraph{Proof:} This proposition can be easily proved using product rule on $\overline{\boldsymbol{v}}=\boldsymbol{Q}^T.\boldsymbol{v}$ and $\overline{\boldsymbol{A}}=\boldsymbol{Q}^T.\boldsymbol{A}.\boldsymbol{Q}$ and substituting for $\partial_{\xi}\boldsymbol{Q}^T=-\boldsymbol{Q}^T.\hat{\boldsymbol{\kappa}}$. The result obtained after such computations (say for $n=1,2,3$) when compared with the results obtained \textit{Proposition 1} and \textit{Corollary 1}, proves the intended result. $\square$
\subsection{Variation and linearization of rotation tensor}
To obtain the virtual rotation tensor field, we superimpose an admissible variation field  $\delta\boldsymbol{Q}$ to the rotation field $\boldsymbol{Q}$. The varied configuration is then defined by $\boldsymbol{Q}_{\epsilon}$ such that for $\epsilon>0$, we have
\begin{equation} \label{varied_Q}
\begin{gathered}
\boldsymbol{Q}_{\epsilon}=\boldsymbol{Q}(\boldsymbol{\uptheta}+\epsilon\delta\boldsymbol{\uptheta})=\boldsymbol{Q}(\epsilon\delta\boldsymbol{\alpha}).\boldsymbol{Q}(\boldsymbol{\uptheta})\\
\delta\boldsymbol Q=\partial_{\epsilon}\boldsymbol{Q}_{\epsilon}|_{\epsilon=0}.
\end{gathered}
\end{equation}
The fact that $SO(3)$ is a non-linear manifold makes it difficult to geometrically understand and obtain the variation of rotation tensor. We also note that it is advantageous to express the virtual rotation tensor by means of\textit{ virtual rotation vector in current state} $\delta\boldsymbol{\alpha}$ (that is saying $\delta\hat{\boldsymbol{\alpha}}.\boldsymbol{Q}\in T_{\boldsymbol{Q}}SO(3)$) contrary to the \textit{variation of total rotation vector} $\delta\boldsymbol{\uptheta}$ ($\delta\hat{\boldsymbol{\uptheta}}\in so(3)$). The varied director field is then given by
\begin{equation} \label{varied_d}
\boldsymbol{d}_{i\epsilon}=\boldsymbol{Q}_{\epsilon}.\boldsymbol{E}_i=\boldsymbol{Q}(\epsilon\delta\boldsymbol{\alpha}).\boldsymbol{d}_i
\end{equation}
The rotation tensor $\boldsymbol{Q}_{\epsilon}=\boldsymbol{Q}(\boldsymbol\uptheta+\epsilon\delta\boldsymbol{\uptheta})$ transforms the vector $\boldsymbol{E}_i$ to $\boldsymbol{d}_{i\epsilon}$ in a single step, whereas, the tensor $\boldsymbol{Q}_{\epsilon}=\boldsymbol{Q}(\epsilon\delta\boldsymbol{\alpha}).\boldsymbol{Q}(\boldsymbol{\uptheta})$ performs the same transformation in two steps: $\boldsymbol{E}_i\xrightarrow{\boldsymbol{Q}(\boldsymbol{\uptheta})}\boldsymbol{d}_i\xrightarrow{\boldsymbol{Q}(\epsilon\delta\boldsymbol{\alpha})}\boldsymbol{d}_{i\epsilon}$. From Eq. \eqref{varied_Q}, we arrive at the expression of varied rotation tensor and director field:
\begin{subequations} 
	\begin{gather}
	\delta\boldsymbol{Q}=\partial_{\epsilon}\left(\text{exp}(\epsilon\delta\hat{\boldsymbol\alpha}).\text{exp}(\hat{\boldsymbol\uptheta})\right)|_{\epsilon=0}=\left(\delta\hat{\boldsymbol{\alpha}}.\text{exp}(\epsilon\delta\hat{\boldsymbol\alpha}).\text{exp}(\hat{\boldsymbol\uptheta})\right)|_{\epsilon=0}=\delta\hat{\boldsymbol{\alpha}}.\boldsymbol{Q}(\boldsymbol{\uptheta});\label{variation_Q}\\
	\delta\boldsymbol{d}_i=\delta\boldsymbol{Q}.\boldsymbol{E}_i=\delta\hat{\boldsymbol{\alpha}}.\boldsymbol{d}_i.\label{variation_d}.
	\end{gather}
\end{subequations}
Here, $\delta\hat{\boldsymbol{\alpha}}$ represents the anti-symmetric matrix associated with the vector  $\delta\boldsymbol{\alpha}$. We define the material form of incremental rotation $\delta\hat{\overline{\boldsymbol{\alpha}}}$ (with $\delta{\overline{\boldsymbol{\alpha}}}$ being the associated axial vector) as
\begin{equation} \label{alpha_bar_hat}
\delta\hat{\overline{\boldsymbol{\alpha}}}=\boldsymbol{Q}^T.\delta\hat{\boldsymbol{\alpha}}=\boldsymbol{Q}^T.\delta\boldsymbol{Q};\text{ }\delta{\overline{\boldsymbol{\alpha}}}=\boldsymbol{Q}^T.\delta\boldsymbol{\alpha}.
\end{equation}
Like the variation, we define the linearized part of rotation tensor $\boldsymbol{Q}$, linearized at $\text{exp}(\hat{\boldsymbol{\uptheta}})$ in the direction of $\Delta\hat{\boldsymbol{\alpha}}.\boldsymbol{Q}\in T_{\boldsymbol{Q}}SO(3)$ as
\begin{equation} \label{l3}
\Delta\boldsymbol{Q}=\partial_{\epsilon}\boldsymbol{Q}_{\epsilon}|_{\epsilon=0} \text{ with } 	\boldsymbol{Q}_{\epsilon}=\boldsymbol{Q}(\epsilon\Delta\boldsymbol{\alpha}).\boldsymbol{Q}(\boldsymbol{\uptheta}).
\end{equation}
 
It follows from the discussion above that $\partial_{\xi}\boldsymbol{Q},\delta\boldsymbol{Q},\Delta\boldsymbol{Q}\in T_{\boldsymbol{Q}}SO(3)$, and $\delta\hat{\overline{\boldsymbol{\alpha}}},\delta\hat{\boldsymbol{\uptheta}},\Delta\hat{\boldsymbol{\uptheta}}\in so(3)$. Like the relationship between $\boldsymbol{\kappa}$ and $\partial_{\xi}\boldsymbol{\uptheta}$, we arrive at the relation between $\delta\boldsymbol{\alpha}$ (or $\Delta\boldsymbol\alpha$) and $\delta\boldsymbol{\uptheta}$ (or $\Delta\boldsymbol\uptheta$). We utilize the results in Eq. \eqref{varied_Q} and obtain
\begin{equation} \label{alpha_uptheta1}
\partial_{\epsilon}\text{exp}(\epsilon\delta\hat{\boldsymbol\alpha})|_{\epsilon=0}=\partial_{\epsilon}\left(\text{exp}(\hat{\boldsymbol\uptheta}+\epsilon\delta\hat{\boldsymbol\uptheta}).\text{exp}(-\hat{\boldsymbol\uptheta})\right)|_{\epsilon=0}.
\end{equation}
Simplifying Eq. \eqref{alpha_uptheta1}, we get
\begin{equation}\label{alpha_uptheta2}
\delta{\boldsymbol\alpha}=\boldsymbol{T}_{\boldsymbol{\uptheta}}.\delta{\boldsymbol\uptheta};\text{ }\delta{\boldsymbol\uptheta}=\boldsymbol{T}_{\boldsymbol{\uptheta}}^{-1}.\delta{\boldsymbol\alpha}.
\end{equation}
\textit{Proposition 3} also holds for the variation and mix of derivatives and variation, for example:  $\delta^n\overline{\boldsymbol{v}}=\boldsymbol{Q}^T.\tilde\delta^n\boldsymbol{v}$  and $\delta^k(\partial_{\xi}^n\overline{\boldsymbol{v}})=\boldsymbol{Q}^T.\tilde\delta^k(\tilde\partial_{\xi}^n\boldsymbol{v})$, where $\tilde\delta^k=(\delta-\hat\delta)^k$ such that $\hat\delta=\delta\hat{\boldsymbol{\alpha}}$.

\begin{figure}[htb]\centering 
	\includegraphics[width=\textwidth]{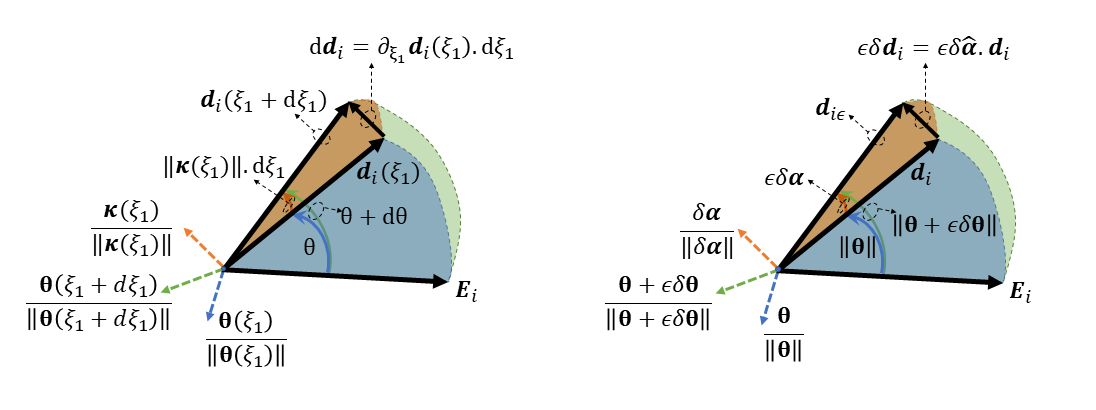}
	\caption{Physical interpretation of curvature $\boldsymbol\kappa$ (left figure) and variation of rotation vector $\delta\boldsymbol\alpha$ (right figure) resulting in infinitesimal rotation}
	\label{fig:Fig1}
\end{figure}
 In figure \ref{fig:Fig1}, we originate three vectors (the reference vector $\boldsymbol{E}_i$, the vector $\boldsymbol{d}_i(\xi)$ obtained by finite rotation of $\boldsymbol{E}_i$, and the vector $\boldsymbol{d}_i(\xi+\text{d}\xi)$) at a point to illustrates the concept of curvature and the incremental (or variation) current rotation vector.
\section{On derivatives}
\label{section3}
\subsection{Useful results on derivatives of Lie-bracket and higher-order product rule}
\paragraph{Proposition 4:} For any $\hat{\boldsymbol{a}},\hat{\boldsymbol{b}}\in so(3)$ with corresponding axial vectors $\boldsymbol{a},\boldsymbol{b}\in \mathbb{R}^3$ respectively, the following formula for derivatives hold:
\begin{subequations}\label{lie_brackets_der}
	\begin{gather}
	\partial^{n}_{\xi}\left[\hat{\boldsymbol{a}},\hat{\boldsymbol{b}}\right]=\sum_{i=0}^n C_i^n\left[\partial_{\xi}^{(n-i)}\hat{\boldsymbol{a}},\partial_{\xi}^{(i)}\hat{\boldsymbol{b}}\right]=\sum_{i=0}^n C_i^n\left[\partial_{\xi}^{(i)}\hat{\boldsymbol{a}},\partial_{\xi}^{(n-i)}\hat{\boldsymbol{b}}\right];\label{lie_brackets_der1}\\
	\partial^{n}_{\xi}({\boldsymbol{a}}\times{\boldsymbol{b}})=\sum_{i=0}^n C_i^n(\partial_{\xi}^{(n-i)}{\boldsymbol{a}}\times\partial_{\xi}^{(i)}{\boldsymbol{b}})=\sum_{i=0}^n C_i^n(\partial_{\xi}^{(i)}{\boldsymbol{a}}\times\partial_{\xi}^{(n-i)}{\boldsymbol{b}});\label{lie_brackets_der2}
	\end{gather}
\end{subequations}
\paragraph{Proof:} Using the definition of Lie-bracket in Eq. \eqref{lie_brackets1}, we have
\begin{equation}\label{lie_brackets3}\partial_{\xi}\left[\hat{\boldsymbol{a}},\hat{\boldsymbol{b}}\right]=(\partial_{\xi}\hat{\boldsymbol{a}}.\hat{\boldsymbol{b}}+\hat{\boldsymbol{a}}.\partial_{\xi}\hat{\boldsymbol{b}})-(\partial_{\xi}\hat{\boldsymbol{b}}.\hat{\boldsymbol{a}}+\hat{\boldsymbol{b}}.\partial_{\xi}\hat{\boldsymbol{a}})=\left[\partial_{\xi}\hat{\boldsymbol{a}},\hat{\boldsymbol{b}}\right]+\left[\hat{\boldsymbol{a}},\partial_{\xi}\hat{\boldsymbol{b}}\right].
\end{equation} 

Higher-order derivatives of the Lie-bracket derived using the above result yields an expression given by Eq. \eqref{lie_brackets_der1}. Using the definition of axial vector corresponding to the Lie-bracket in Eq. \eqref{lie_brackets2}, Eq. \eqref{lie_brackets_der2} follows from Eq. \eqref{lie_brackets_der1}. The first and second equality in \eqref{lie_brackets_der1} and \eqref{lie_brackets_der2} holds by virtue of result \eqref{binomial_prop3} in \textit{Theorem 0}.
This completes the proof. $\square$
\paragraph{Corollary 2:} For scalars functions $f(\xi),g(\xi)$; vectors $\boldsymbol{u}(\xi),\boldsymbol{v}(\xi)$; and a second order tensors $\boldsymbol{A}(\xi), \boldsymbol{B}(\xi)$, the $n^{\text{th}}$ order product rule is given by the following:
\begin{subequations}\label{general_der}
	\begin{align}
	&\partial^{n}_{\xi}(fg)=\sum_{i=0}^n C_i^n\partial_{\xi}^{(n-i)}f.\partial_{\xi}^{(i)}g=\sum_{i=0}^n C_i^n\partial_{\xi}^{(i)}f.\partial_{\xi}^{(n-i)}g;\label{general_der1}\\
	&\partial^{n}_{\xi}(f\boldsymbol{u})=\sum_{i=0}^n C_i^n\partial_{\xi}^{(n-i)}f.\partial_{\xi}^{(i)}\boldsymbol{u}=\sum_{i=0}^n C_i^n\partial_{\xi}^{(i)}f.\partial_{\xi}^{(n-i)}\boldsymbol{u};\label{general_der2}\\
	&\partial^{n}_{\xi}(f\boldsymbol{A})=\sum_{i=0}^n C_i^n\partial_{\xi}^{(n-i)}f.\partial_{\xi}^{(i)}\boldsymbol{A}=\sum_{i=0}^n C_i^n\partial_{\xi}^{(i)}f.\partial_{\xi}^{(n-i)}\boldsymbol{A};\label{general_der3}\\
	&\partial^{n}_{\xi}(\boldsymbol{v}\cdot\boldsymbol{u})=\sum_{i=0}^n C_i^n\partial_{\xi}^{(n-i)}\boldsymbol{v}\cdot\partial_{\xi}^{(i)}\boldsymbol{u}=\sum_{i=0}^n C_i^n\partial_{\xi}^{(i)}\boldsymbol{v}\cdot\partial_{\xi}^{(n-i)}\boldsymbol{u};\label{general_der4}\\
	&\partial^{n}_{\xi}(\boldsymbol{A}.\boldsymbol{u})=\sum_{i=0}^n C_i^n\partial_{\xi}^{(n-i)}\boldsymbol{A}.\partial_{\xi}^{(i)}\boldsymbol{u}=\sum_{i=0}^n C_i^n\partial_{\xi}^{(i)}\boldsymbol{A}.\partial_{\xi}^{(n-i)}\boldsymbol{u}.\label{general_der5}\\
	&\partial^{n}_{\xi}(\boldsymbol{A}.\boldsymbol{B})=\sum_{i=0}^n C_i^n\partial_{\xi}^{(n-i)}\boldsymbol{A}.\partial_{\xi}^{(i)}\boldsymbol{B}=\sum_{i=0}^n C_i^n\partial_{\xi}^{(i)}\boldsymbol{A}.\partial_{\xi}^{(n-i)}\boldsymbol{B}.\label{general_der6}
	\end{align}
\end{subequations}
\paragraph{Proof:} The result above follows directly from proposition 4. This is because, from Eq. \eqref{lie_brackets3}, we see that Lie-brackets follow the chain rule just like the product of two scalars or a dot product except for the fact that Lie-brackets are non-commutative (which is a stronger condition). $\square$ 
\paragraph{Proposition 5:} For any $\hat{\boldsymbol{a}}(\xi)\in so(3)$, the following holds: 
\begin{equation}\label{lie_brackets_der3}
	\partial^{m}_{\xi}\left[\hat{\boldsymbol{a}},\partial_{\xi}\hat{\boldsymbol{a}}\right]=\sum_{j=0}^{j_{m}^{\text{max}}} b_{mj}\left[\partial_{\xi}^{(j)}\hat{\boldsymbol a},\partial_{\xi}^{(m-j+1)}\hat{\boldsymbol a}\right].
\end{equation}
where $m,j,j^{\text{max}}_{m},b_{mj}\in\mathbb{Z}^+$, such that the coefficient $j^{\text{max}}_{m}$, and $b_{mj}$ are given as
\begin{subequations}\label{lie_brackets_der4}
\begin{align}
&j^{\text{max}}_{m}=\begin{cases}\text{floor}\left(\frac{m+1}{2}\right), & \text{ if } \frac{m+1}{2}\notin \mathbb{Z}^+;\\
\frac{m+1}{2}-1,& \text{ if } \frac{m+1}{2}\in \mathbb{Z}^+.
\end{cases}\label{lie_brackets_der4a}\\
&b_{mj}=C_{j}^m-C_{(m-j+1)}^m=C^m_j\left(\frac{m-2j+1}{m-j+1}\right)\label{lie_brackets_der4b}.
\end{align}
\end{subequations}
\paragraph{Proof:}  Using the result \eqref{lie_brackets_der1} of \textit{Proposition 4}, we get,
\begin{equation}\label{lie_brackets_der5}
\begin{aligned}
\partial_{\xi}^{m}[\hat{\boldsymbol{a}},\partial_{\xi}\hat{\boldsymbol{a}}]&=\sum_{j=0}^{m}C_{j}^m\left[\partial_{\xi}^{(j)}\hat{\boldsymbol{a}},\partial_{\xi}^{(m-j+1)}\hat{\boldsymbol{a}}\right]\\
&=\underbrace{C^m_0\left[\hat{\boldsymbol{a}},\partial_{\xi}^{(m+1)}\hat{\boldsymbol{a}}\right]}_{\text{Term 0}}+\underbrace{\sum_{j=1}^{m}C_{j}^m\left[\partial_{\xi}^{(j)}\hat{\boldsymbol{a}},\partial_{\xi}^{(m-j+1)}\hat{\boldsymbol{a}}\right]}_{\text{Term 1}}.
\end{aligned}
\end{equation}
The terms in the expansion of \textit{Term 1} with $j=\frac{m+1}{2}$ vanishes (refer to Eq. \eqref{lie_brackets_properties1}). Keeping that in mind, we note that the expansion of \textit{Term 1} can be written in two possible ways: the first possibility is when $j>\frac{m+1}{2}$, and the second option is considering $j<\frac{m+1}{2}$. In the first option, the coefficient of all the terms in the sum will be negative, whereas, for second case, the coefficients will be positive. This owes to the anti-commutative property of Lie-brackets mentioned in Eq. \eqref{lie_brackets_properties2}. We consider the second case in our derivation. 

We can further simplify \textit{Term 1}. The total number of terms present in the expanded form of \textit{Term 1} is less than $m$. This is because the terms with interchanged order of derivatives in Lie-bracket can be reduced to one term. For instance: $$c_1\left[\partial_{\xi}^{x}\hat{\boldsymbol{a}},\partial_{\xi}^{y}\hat{\boldsymbol{a}}\right]+c_2\left[\partial_{\xi}^{y}\hat{\boldsymbol{a}},\partial_{\xi}^{x}\hat{\boldsymbol{a}}\right]=(c_1-c_2)\left[\partial_{\xi}^{x}\hat{\boldsymbol{a}},\partial_{\xi}^{y}\hat{\boldsymbol{a}}\right].$$ Thus, the maximum value of $j$ is restricted by the fact that $j<\frac{m+1}{2}$ and $j\in\mathbb{Z}^+-\{0\}$. Combining these two constraints yields $\text{max}(j)=j^{\text{max}}_m$ given by Eq. \eqref{lie_brackets_der4a}. However, such a reduction or simplification of \textit{Term 1} changes the coefficient by which each term in the sum is weighed. The discussion presented so far may be demonstrated, for example, for $m=4$ as: 
\begin{equation}\label{lie_brackets_der6}
\begin{aligned}
\text{Term 1}|_{m=4}&=\sum_{j=1}^{4}C_{j}^4\left[\partial_{\xi}^{(j)}\hat{\boldsymbol{a}},\partial_{\xi}^{(5-j)}\hat{\boldsymbol{a}}\right]\\&=3\left[\partial_{\xi}\hat{\boldsymbol{a}},\partial_{\xi}^4\hat{\boldsymbol{a}}\right]+2\left[\partial_{\xi}^2\hat{\boldsymbol{a}},\partial_{\xi}^3\hat{\boldsymbol{a}}\right]=-3\left[\partial_{\xi}^4\hat{\boldsymbol{a}},\partial_{\xi}\hat{\boldsymbol{a}}\right]-2\left[\partial_{\xi}^3\hat{\boldsymbol{a}},\partial_{\xi}^2\hat{\boldsymbol{a}}\right]\\
&=\|C_{1}^4-C_{4}^4\|\left[\partial_{\xi}\hat{\boldsymbol{a}},\partial_{\xi}^4\hat{\boldsymbol{a}}\right]+\|C_{2}^4-C_{3}^4\|\left[\partial_{\xi}^2\hat{\boldsymbol{a}},\partial_{\xi}^3\hat{\boldsymbol{a}}\right]\\
&=\sum_{j=1}^{j_{4}^{max}=2}(C_{j}^4-C_{(4-j+1)}^4)\left[\partial_{\xi}^{(j)}\hat{\boldsymbol{a}},\partial_{\xi}^{(5-j)}\hat{\boldsymbol{a}}\right]
\end{aligned}
\end{equation}
For a general case, if $j_{m}^{\text{max}}\geq 1$, \textit{Term 1} can be written as:
\begin{equation}\label{lie_brackets_der7}
\text{Term 1}=\sum_{j=1}^{j_{m}^{\text{max}}} b_{mj}\left[\partial_{\xi}^{(j)}\hat{\boldsymbol a},\partial_{\xi}^{(m-j+1)}\hat{\boldsymbol a}\right].
\end{equation}
The modified coefficient $b_{mj}=(C_{j}^m-C_{(m-j+1)}^m)$ is defined in Eq. \eqref{lie_brackets_der4b}. From the second equality in Eq. \eqref{lie_brackets_der4b}, we also note that $b_{m0}=C^m_0=1$. Therefore, combining \textit{Term 0} and \textit{Term 1} proves the proposition. Table \ref{table2} gives the value of coefficient $j^{\text{max}}_{m}$ for $1\leq m\leq6$. $\square$
\begin{table}[htb]
	\centering
	\begin{tabular}{|l|lllllll|}
		\hline	&$i=0$ & $i=1$ & $i=2$ & $i=3$ & $i=4$ & $i=5$ & $i=6$ \\\hline
		$n=0$ & 0 & - & - & - & - & - & - \\\hline
		$n=1$ & 0 & 0 & - & - & - & - & - \\\hline
		$n=2$ & 1 & 0 & 0 & - & - & - & - \\\hline
		$n=3$ & 1 & 1 & 0 & 0 & - & - & - \\\hline
		$n=4$ & 2 & 1 & 1 & 0 & 0 & - & - \\\hline
		$n=5$ & 2 & 2 & 1 & 1 & 0 & 0 & - \\\hline
		$n=6$ & 3 & 2 & 2 & 1 & 1 & 0 & 0\\\hline
	\end{tabular}
	\caption{$j^{\text{max}}_{m}$ for $0\leq n \leq 6$}
	\label{table2}
\end{table}
\FloatBarrier
\subsection{Derivatives of curvature tensor}
The derivative of the curvature tensor may be obtained using Eq. \eqref{curv_tensor_uptheta} deploying a straightforward application of the chain rule. However, deriving the expression of higher-order derivatives using Eq. \eqref{curv_tensor_uptheta} is cumbersome because of involvement of trigonometric functions. Instead, we realize that the reparametrization of the rotation tensor by the \textit{Gibbs vector} (the components of which are called as \textit{Gibbs} or \textit{Rodriguez parameters} in the literature) yields the formula of curvature tensor that is beneficial in obtaining the derivative of curvature tensor in the form of a single summation-formula. Consider a rotation tensor $\boldsymbol{Q}(\boldsymbol\uptheta)=\text{exp}(\hat{\boldsymbol\uptheta})\in SO(3)$. We define the Gibbs vector $\boldsymbol{\phi}$ and the associated quantities as:
\begin{equation}\label{l34}
\begin{gathered}
\hat{\boldsymbol{\phi}}=\frac{\tan\left(\frac{\uptheta}{2}\right)}{\uptheta}\hat{\boldsymbol\uptheta};\quad \boldsymbol{\phi}=\frac{\tan\left(\frac{\uptheta}{2}\right)}{\uptheta}{\boldsymbol\uptheta};\quad \phi=\|\boldsymbol{\phi}\|=\tan\left(\frac{\uptheta}{2}\right);\quad \overline\phi=2\cos^2\left(\frac{\uptheta}{2}\right)=\frac{2}{\phi^2+1}.
\end{gathered}
\end{equation}
The result defined in Eq. \ref{exp_map1} may be manipulated using the definition given above as:
\begin{equation}\label{l35}
\begin{aligned}
&\boldsymbol{Q}(\hat{\boldsymbol{\phi}})=\boldsymbol{I}_3+2\cos^2\left(\frac{\uptheta}{2}\right).\left(\frac{1}{\uptheta}\tan\left(\frac{\uptheta}{2}\right)\hat{\boldsymbol\uptheta}+\left(\frac{1}{\uptheta}\tan\left(\frac{\uptheta}{2}\right)\right)^2\hat{\boldsymbol\uptheta}^2\right)=\boldsymbol{I}_3+\overline\phi(\hat{\boldsymbol{\phi}}+\hat{\boldsymbol{\phi}}^2);\\
&\boldsymbol{Q}(\hat{\boldsymbol{\phi}})^T=\boldsymbol{Q}(-\boldsymbol{\uptheta})=\boldsymbol{I}_3+\overline\phi(-\hat{\boldsymbol{\phi}}+\hat{\boldsymbol{\phi}}^2).
\end{aligned}
\end{equation}

\paragraph{Proposition 6:} $\hat{\boldsymbol{\phi}}\in so(3)$ and $n\in\mathbb{Z}^+-\{0\}$, the following formulae hold true:\\
\begin{subequations} 
	\begin{align}
	&\hat{\boldsymbol{\phi}}^{2n-1}=(-1)^{n-1}\phi^{2(n-1)}\hat{\boldsymbol{\phi}};\quad\hat{\boldsymbol{\phi}}^{2n}=(-1)^{n-1}\phi^{2(n-1)}\hat{\boldsymbol{\phi}}^2;\label{idset41}\\
	&\hat{\boldsymbol{\phi}}.\partial_{\xi}\hat{\boldsymbol{\phi}}.\hat{\boldsymbol{\phi}}=-(\boldsymbol{\phi}\cdot\partial_{\xi}\boldsymbol{\phi})\hat{\boldsymbol{\phi}};\quad \hat{\boldsymbol{\phi}}.\partial_{\xi}\hat{\boldsymbol{\phi}}.\hat{\boldsymbol{\phi}}^2=-(\boldsymbol{\phi}\cdot\partial_{\xi}\boldsymbol{\phi})\hat{\boldsymbol{\phi}}^2.\label{idset42}
	\end{align}
\end{subequations}
\paragraph{Proof:} Refer to Eq. [30] of Argyris \cite{argyris1982excursion} for identity \eqref{idset41} that describes the recurrence formula for the power of anti-symmetric matrix. The identity \eqref{idset42} can be proved by considering the action of tensor on left hand side of equation on to a vector, say $\boldsymbol{v}\in\mathbb{R}^3$ and using the vector triple product identity, such that
\begin{equation}\label{l37}
\begin{aligned}
(\hat{\boldsymbol{\phi}}.\partial_{\xi}\hat{\boldsymbol{\phi}}.\hat{\boldsymbol{\phi}}).\boldsymbol{v}&=\boldsymbol{\phi}\times(\partial_{\xi}{\boldsymbol{\phi}}\times(\hat{\boldsymbol{\phi}}.\boldsymbol{v}))={({\boldsymbol{\phi}}\cdot(\hat{\boldsymbol{\phi}}.\boldsymbol{v}))\partial_{\xi}\boldsymbol{\phi}}-(\boldsymbol{\phi}\cdot\partial_{\xi}\boldsymbol{\phi})\hat{\boldsymbol{\phi}}.\boldsymbol{v}.
\end{aligned}
\end{equation}
Noting that $({\boldsymbol{\phi}}\cdot(\hat{\boldsymbol{\phi}}.\boldsymbol{v}))=\boldsymbol{\phi}\cdot(\boldsymbol{\phi}\times\boldsymbol{v})=\boldsymbol{0}$, we prove the first part of identity \eqref{idset42}. Along the similar lines, the second part can be proven. $\square$

The curvature tensor can then be obtained using Eq. \eqref{l35} and \textit{proposition 6} as:
\begin{equation}\label{l36}
\hat{\boldsymbol{\kappa}}=\partial_{\xi}\boldsymbol{Q}.\boldsymbol{Q}^T=\overline{\phi}\left(\partial_{\xi}\hat{\boldsymbol{\phi}}+[\hat{\boldsymbol{\phi}},\partial_{\xi}\hat{\boldsymbol{\phi}}]\right).
\end{equation}
This expression is much simpler than the one presented in Eq. \eqref{curv_tensor_uptheta}. 
\paragraph{Proposition 7:}The following hold:
\begin{subequations}\label{l38}
\begin{align}
\partial_{\xi}^n\hat{\boldsymbol\kappa}&=\sum_{i=0}^{n}C_{i}^n\partial_{\xi}^{(i)}\overline{\phi}\left(\partial_{\xi}^{(n-i+1)}\hat{\boldsymbol{\phi}}+\sum_{j=0}^{j^{\text{max}}_{(n-i)}}b_{(n-i)j}\left[\partial_{\xi}^{(j)}\hat{\boldsymbol{\phi}},\partial_{\xi}^{(n-i+1-j)}\hat{\boldsymbol{\phi}}\right]\right);\label{l38a}\\
\partial_{\xi}^n{\boldsymbol\kappa}&=\sum_{i=0}^{i=n}C_{i}^n\partial_{\xi}^{(i)}\overline{\phi}\left(\partial_{\xi}^{(n-i+1)}{\boldsymbol{\phi}}+\sum_{j=0}^{j^{\text{max}}_{(n-i)}}b_{(n-i)j}\left(\partial_{\xi}^{(j)}{\boldsymbol{\phi}}\times\partial_{\xi}^{(n-i+1-j)}{\boldsymbol{\phi}}\right)\right).\label{l38b}
\end{align}
\end{subequations}
where, $n,i,j,j^{\text{max}}_{(n-i)},C_{i}^n,b_{(n-i)j}\in\mathbb{Z}^+$. Replacing $m\longrightarrow (n-i)$ in \eqref{lie_brackets_der4a} and \eqref{lie_brackets_der4b} yields $j^{\text{max}}_{(n-i)}$, and $b_{(n-i)j}$.
\paragraph{Proof:}We utilize Eq.\eqref{general_der3} of \textit{Corollary 2} and the expression of curvature tensor in Eq. \eqref{l36} and obtain
\begin{equation}\label{l39a}
\begin{aligned}
\partial_{\xi}^n\hat{\boldsymbol\kappa}=&\sum_{i=0}^n C_{i}^n\partial_{\xi}^{(i)}\overline{\phi}.\left(\partial_{\xi}^{(n-i)}\left(\partial_{\xi}\hat{\boldsymbol{\phi}}+\left[\hat{\boldsymbol{\phi}},\partial_{\xi}\hat{\boldsymbol{\phi}}\right]\right)\right)\\
=&\sum_{i=0}^n C_{i}^n\partial_{\xi}^{(i)}\overline{\phi}.\left(\partial_{\xi}^{(n-i+1)}\hat{\boldsymbol{\phi}}+\partial_{\xi}^{(n-i)}\left[\hat{\boldsymbol{\phi}},\partial_{\xi}\hat{\boldsymbol{\phi}}\right]\right)
\end{aligned}
\end{equation}
Using \textit{Proposition 5} and replacing $m\longrightarrow (n-i)$, we get,
\begin{equation}\label{l39b}
\partial^{(n-i)}_{\xi}\left[\hat{\boldsymbol{\phi}},\partial_{\xi}\hat{\boldsymbol{\phi}}\right]=\sum_{j=0}^{j_{(n-i)}^{\text{max}}} b_{(n-i)j}\left[\partial_{\xi}^{(j)}\hat{\boldsymbol\phi},\partial_{\xi}^{((n-i)-j+1)}\hat{\boldsymbol\phi}\right].
\end{equation}
The equation above when substituted into Eq. \eqref{l39a} proves the result \eqref{l38a}. The axial vector corresponding to $\partial_{\xi}^n\hat{\boldsymbol\kappa}$, given by Eq. \eqref{l38b}, is obtained using the formula \eqref{lie_brackets2}.  $\square$

Appendix \ref{A61} gives a MATLAB code to obtain spatial curvature vector and its derivatives in accordance with \textit{Proposition 7}.

\paragraph{Note:} Define $\frac{\boldsymbol{\uptheta}}{\uptheta}=\boldsymbol{e}$. We can use Eq. \eqref{l35} and \eqref{l38} to obtain $\boldsymbol{Q}$ and $\partial_{\xi}^n\hat{\boldsymbol{\kappa}}$ for small rotations (when $\|\boldsymbol\uptheta\|\to0$), by setting $\partial_{\xi}^n\phi=\lim_{\uptheta\to 0}\partial_{\xi}^n\tan\left(\frac{\uptheta}{2}\right)$,  $\partial_{\xi}^n\hat{\boldsymbol\phi}=\left(\lim_{\uptheta\to 0}\partial_{\xi}^n\tan\left(\frac{\uptheta}{2}\right)\right).\hat{\boldsymbol{e}}$ and $\partial_{\xi}^n\overline\phi=2\lim_{\uptheta\to 0}\partial_{\xi}^n\cos^2\left(\frac{\uptheta}{2}\right)$. Here, $\boldsymbol{e}$ is the fixed unit vector about which the rotation occurs.
\paragraph{Corollary 3:} The following holds:
\begin{subequations}\label{c2_1}
	\begin{align}
	&\begin{aligned}
	\tilde\partial_{\xi}^n\hat{\boldsymbol{\kappa}}&=\partial_{\xi}^n\hat{\boldsymbol{\kappa}}-(1-\delta_{n0})\sum_{i=1}^{n-1}\partial_{\xi}^{(i-1)}\left[\hat{\boldsymbol{\kappa}},\tilde\partial_{\xi}^{(n-i)}\hat{\boldsymbol \kappa}\right]\\&=\partial_{\xi}^n\hat{\boldsymbol{\kappa}}-(1-\delta_{n0})\sum_{i=1}^{n-1}\sum_{j=0}^{i-1}C_j^{(i-1)}\left[\partial_{\xi}^{(j)}\hat{\boldsymbol{\kappa}},\partial_{\xi}^{(n-i-j)}\tilde\partial_{\xi}^{(n-i)}\hat{\boldsymbol \kappa}\right];
	\end{aligned}\label{c2_1a}\\
	&\begin{aligned}
	\tilde\partial_{\xi}^n{\boldsymbol{\kappa}}=(\partial_{\xi}-\hat\partial)^n\boldsymbol{\kappa}&=\partial_{\xi}^n{\boldsymbol{\kappa}}-(1-\delta_{n0})\sum_{i=1}^{n-1}\partial_{\xi}^{(i-1)}({\boldsymbol{\kappa}}\times\tilde\partial_{\xi}^{(n-i)}{\boldsymbol {\kappa}})\\
	&=\partial_{\xi}^n{\boldsymbol{\kappa}}-(1-\delta_{n0})\sum_{i=1}^{n-1}\sum_{j=0}^{i-1}C^{(i-1)}_{j}(\partial_{\xi}^{(j)}{\boldsymbol{\kappa}}\times\partial_{\xi}^{(n-i-j)}\tilde\partial_{\xi}^{(n-i)}{\boldsymbol {\kappa}}).
	\end{aligned}\label{pc_1b}\\
	&\partial_{\xi}^n\hat{\overline{\boldsymbol{\kappa}}}=\boldsymbol{Q}^T.\tilde\partial_{\xi}^n\hat{\boldsymbol{\kappa}}.\boldsymbol{Q}\label{pc_1c}
	\end{align}
\end{subequations}
\paragraph{Proof:} This corollary follows from the Proposition 1, 2, 3 and 4. We also note that in the sums presented above, $\text{max}(i)=(n-1)$, because $\left[\hat{\boldsymbol{\kappa}},\tilde\partial_{\xi}^{(n-i)}\hat{\boldsymbol{\kappa}}\right]\big|_{(n=i)}=\boldsymbol{0}_3$. $\square$ 

The n$^\text{th}$ order derivative of rotation tensor $\boldsymbol{Q}$ can be derived as a function of Gibbs vector and the associated parameters using Eq. \eqref{l35}. However, computationally, a much simpler approach would be to derive a recurrence formula for $\partial_{\xi}^n\boldsymbol{Q}$ using the fact that $\partial_{\xi}\boldsymbol{Q}=\hat{\boldsymbol\kappa}.\boldsymbol{Q}$ and \textit{Proposition 7}. The recurrence formula for $\partial_{\xi}^n\boldsymbol{Q}$ yields the formula to obtain n$^\text{th}$ order derivative of director vectors $\boldsymbol{d}_m(\xi)$ with $m\in\{1,2,3\}$.
\paragraph{Proposition 8:} For $n\geq0$, the following hold:
\begin{subequations}\label{c2_2}
	\begin{align}
	&\partial_{\xi}^n\boldsymbol{Q}=\delta_{n0}\boldsymbol{Q}+(1-\delta_{n0})\sum_{i=0}^{n-1}C_i^{(n-1)}\partial_{\xi}^i\hat{\boldsymbol\kappa}.\partial_{\xi}^{(n-1-i)}\boldsymbol{Q};\label{c2_2a}\\
	&\partial_{\xi}^n\boldsymbol{d}_m=\delta_{n0}\boldsymbol{d}_m+(1-\delta_{n0})\sum_{i=0}^{n-1}C_i^{(n-1)}\partial_{\xi}^i\hat{\boldsymbol\kappa}.\partial_{\xi}^{(n-1-i)}\boldsymbol{d}_m.\label{c2_2b}
	\end{align}
\end{subequations}
\paragraph{Proof:}  From the definition of curvature tensor, we have $\partial_{\xi}\boldsymbol{Q}=\hat{\boldsymbol\kappa}.\boldsymbol{Q}$. Therefore, for $n>0$, we have $\partial_{\xi}^n\boldsymbol{Q}=\partial_{\xi}^{(n-1)}\left(\hat{\boldsymbol\kappa}.\boldsymbol{Q}\right)$, which when simplified using Eq.  \eqref{general_der6} yields the result \eqref{c2_2a}. The result \eqref{c2_2b} follows from Eq. \eqref{c2_2a} and the fact that $\partial_{\xi}^n\boldsymbol{d}_m=\partial_{\xi}^n\left(\boldsymbol{Q}.\boldsymbol{E}_m=\partial_{\xi}^n\boldsymbol{Q}.\boldsymbol{E}_m\right)$. $\square$ 
\paragraph{Corollary 4:} Alternate to \textit{Proposition 1} and \textit{Corollary 3}, the quantities $\partial_{\xi}^n\overline{\boldsymbol{\kappa}}$ and $\tilde\partial_{\xi}^n\overline{\boldsymbol{\kappa}}$ can be obtained using the relationship $\overline{\boldsymbol{\kappa}}=\boldsymbol{Q}^T.\boldsymbol{\kappa}$ as:
\begin{subequations}\label{c2_3}
	\begin{align}
	&\partial_{\xi}^n\overline{\boldsymbol{\kappa}}=\sum_{i=0}^{n}C_i^{n}\partial_{\xi}^i\boldsymbol{Q}^T.\partial_{\xi}^{(n-i)}\boldsymbol{\kappa} \qquad \text{for }n\geq 0;\label{c2_3a}\\
	&\tilde\partial_{\xi}^n{\boldsymbol{\kappa}}=\boldsymbol{Q}.\partial_{\xi}^n\overline{\boldsymbol{\kappa}}\qquad\text{for }n> 0.\label{c2_3b}
	\end{align}
\end{subequations}
We use \textit{Corollary 4} to develop a MATLAB code in Appendix \ref{A62} that obtains co-rotational derivatives, material curvature and its derivatives, provided $\partial_{\xi}^{n}\boldsymbol{\kappa}$ and $\partial_{\xi}^n\boldsymbol{Q}$ are known.

\section{Updating the curvature and its derivatives}
\label{section4}
In this section, we shall address the situation where the space curve is evolving with time in steps, such that the transformed curve is also parameterized spatially by $\xi$. At time $t$, let the \textit{initial} rotation tensor field be $\boldsymbol{Q}(\xi,t)=\boldsymbol{Q}_\text{i}(\xi)\in SO(3)$ and in the next time step $(t+1)$, the \textit{updated} (or\textit{ final}) rotation tensor field is $\boldsymbol{Q}(\xi,t+1)=\boldsymbol{Q}_\text{f}(\xi)\in SO(3)$. We assume \textit{Eulerian updating} of rotation tensor field, i.e. the change in rotation tensor field from discrete time step $t$ to $(t+1)$ is given by an \textit{incremental current rotation vector field} $\Delta{\boldsymbol\alpha}$, such that $\Delta\hat{\boldsymbol\alpha}.\boldsymbol{Q}_\text{i}\in T_{\boldsymbol{Q}_\text{i}}SO(3)$. We are given the derivative fields $\partial_{\xi}^n\Delta\boldsymbol{\alpha}$ (or $\partial_{\xi}^n\Delta\hat{\boldsymbol{\alpha}}$) and $\partial_{\xi}^n\boldsymbol{Q}_\text{i}$ up to order $n$ (or equivalently, $\partial_{\xi}^{(n-1)}\hat{\boldsymbol\kappa}_\text{i}$, where $\hat{\boldsymbol\kappa}_\text{i}=\partial_{\xi}\boldsymbol{Q}_\text{i}.\boldsymbol{Q}_\text{i}^T$). The question posed is thus: ``How do we obtain the updated curvature tensor, its spatial, material and co-rotational derivatives up to order $(n-1)$?'' It is clear from \eqref{l38} that the $n^{\text{th}}$-order derivative of the curvature tensor requires up to the $(n+1)^{\text{th}}$ derivative of the corresponding rotation vector. To proceed, we first present the updated rotation tensor as
\begin{equation} \label{l33}
\boldsymbol{Q}_\text{f}=\text{exp}(\Delta\hat{\boldsymbol\alpha}).\boldsymbol{Q}_\text{i}=\boldsymbol{Q}_+.\boldsymbol{Q}_\text{i} \text{ where, } \boldsymbol{Q}_+=\text{exp}(\Delta\hat{\boldsymbol\alpha}).
\end{equation}
We define the curvature corresponding to \textit{incremental current rotation vector} $\Delta\boldsymbol{\alpha}$ and the \textit{transport operator} $\mathbb{T}_{\boldsymbol{Q}}$ as:
\begin{equation}\label{l42}
\begin{gathered}
\hat{\boldsymbol{\kappa}}_{+}=\partial_{\xi}\text{exp}(\Delta\hat{\boldsymbol\alpha}).\text{exp}(-\Delta\hat{\boldsymbol\alpha})=\partial_{\xi}\boldsymbol{Q}_{+}.\boldsymbol{Q}_{+}^T;\\
\mathbb{T}_{\boldsymbol{Q}}[\hat{\boldsymbol{A}}]=\boldsymbol{Q}.\hat{\boldsymbol{A}}.\boldsymbol{Q}^T\in so(3),\quad \forall\quad \boldsymbol{Q}\in SO(3), \hat{\boldsymbol{A}}\in so(3).
\end{gathered}
\end{equation}
We observe that $\mathbb{T}_{\boldsymbol{Q}}[\hat{\overline{\boldsymbol{\kappa}}}]=\hat{\boldsymbol{\kappa}}$ and $\mathbb{T}_{\boldsymbol{Q}^T}[\hat{{\boldsymbol{\kappa}}}]=\hat{\overline{\boldsymbol{\kappa}}}$. 
\paragraph{Proposition 9:} The $n^{\text{th}}$ order derivative of the transport operator $\mathbb{T}_{\boldsymbol{Q}}[\hat{\boldsymbol{A}}]$ is given by
\begin{equation}\label{l42a}
\begin{gathered}
\partial_{\xi}^n\mathbb{T}_{\boldsymbol{Q}}[\hat{\boldsymbol{A}}]=\mathbb{T}_{\boldsymbol{Q}}[\partial_{\xi}^n\hat{\boldsymbol{A}}]+(1-\delta_{n0})\sum_{k=1}^{n}\sum_{i=0}^{n-k}C_{i}^{(n-k)}\left[\partial_{\xi}^{(i)}\hat{\boldsymbol{\kappa}},\partial_{\xi}^{(n-k-i)}\mathbb{T}_{\boldsymbol{Q}}[\partial_{\xi}^{(k-1)}\hat{\boldsymbol{A}}]\right].
\end{gathered}
\end{equation}
\paragraph{Proof:}
Consider
\begin{equation}\label{l42b}
\begin{gathered}
\partial_{\xi}\mathbb{T}_{\boldsymbol{Q}}[\hat{\boldsymbol{A}}]=\boldsymbol{Q}.\partial_{\xi}\hat{\boldsymbol{A}}.\boldsymbol{Q}^T+\hat{\boldsymbol{\kappa}}.\mathbb{T}_{\boldsymbol{Q}}[\hat{\boldsymbol{A}}]-\mathbb{T}_{\boldsymbol{Q}}[\hat{\boldsymbol{A}}].\hat{\boldsymbol{\kappa}}=\mathbb{T}_{\boldsymbol{Q}}[\partial_{\xi}\hat{\boldsymbol{A}}]+[\hat{\boldsymbol{\kappa}},\mathbb{T}_{\boldsymbol{Q}}[\hat{\boldsymbol{A}}]].
\end{gathered}
\end{equation}
Using the above result along with \textit{Proposition 4}, for $n\geq1$, we have
\begin{align*}
\partial_{\xi}^n\mathbb{T}_{\boldsymbol{Q}}[\hat{\boldsymbol{A}}]=&\partial_{\xi}^{(n-1)}.(\partial_{\xi}\mathbb{T}_{\boldsymbol{Q}}[\hat{\boldsymbol{A}}])=\partial_{\xi}^{(n-1)}.\left(\mathbb{T}_{\boldsymbol{Q}}[\partial_{\xi}\hat{\boldsymbol{A}}]+[\hat{\boldsymbol{\kappa}},\mathbb{T}_{\boldsymbol{Q}}[\hat{\boldsymbol{A}}]]\right)\\
=&\partial_{\xi}^{(n-1)}\mathbb{T}_{\boldsymbol{Q}}[\partial_{\xi}\hat{\boldsymbol{A}}]+\sum_{i=0}^{(n-1)}C_{i}^{(n-1)}\left[\partial_{\xi}^{(i)}\hat{\boldsymbol{\kappa}},\partial_{\xi}^{(n-1-i)}\mathbb{T}_{\boldsymbol{Q}}[\hat{\boldsymbol{A}}]\right]\\
=&\partial_{\xi}^{(n-2)}(\partial_{\xi}\mathbb{T}_{\boldsymbol{Q}}[\partial_{\xi}\hat{\boldsymbol{A}}])+\sum_{i=0}^{(n-1)}C_{i}^{(n-1)}\left[\partial_{\xi}^{(i)}\hat{\boldsymbol{\kappa}},\partial_{\xi}^{(n-1-i)}\mathbb{T}_{\boldsymbol{Q}}[\hat{\boldsymbol{A}}]\right]\\
=&\partial_{\xi}^{(n-2)}\mathbb{T}_{\boldsymbol{Q}}[\partial_{\xi}^2\hat{\boldsymbol{A}}]+\sum_{i=0}^{(n-2)}C_{i}^{(n-2)}\left[\partial_{\xi}^{(i)}\hat{\boldsymbol{\kappa}},\partial_{\xi}^{(n-2-i)}\mathbb{T}_{\boldsymbol{Q}}[\partial_{\xi}\hat{\boldsymbol{A}}]\right]\\&+\sum_{i=0}^{(n-1)}C_{i}^{(n-1)}\left[\partial_{\xi}^{(i)}\hat{\boldsymbol{\kappa}},\partial_{\xi}^{(n-1-i)}\mathbb{T}_{\boldsymbol{Q}}[\hat{\boldsymbol{A}}]\right]\\
=&\mathbb{T}_{\boldsymbol{Q}}[\partial_{\xi}^n\hat{\boldsymbol{A}}]+\Bigg(\sum_{i=0}^{(n-1)}C_{i}^{(n-1)}\left[\partial_{\xi}^{(i)}\hat{\boldsymbol{\kappa}},\partial_{\xi}^{(n-1-i)}\mathbb{T}_{\boldsymbol{Q}}[\hat{\boldsymbol{A}}]\right]\\&\sum_{i=0}^{(n-2)}C_{i}^{(n-2)}\left[\partial_{\xi}^{(i)}\hat{\boldsymbol{\kappa}},\partial_{\xi}^{(n-2-i)}\mathbb{T}_{\boldsymbol{Q}}[\partial_{\xi}\hat{\boldsymbol{A}}]\right]+...+\sum_{i=0}^{0}C_{i}^0\left[\partial_{\xi}^{(i)}\hat{\boldsymbol{\kappa}},\partial_{\xi}^{(n-i)}\mathbb{T}_{\boldsymbol{Q}}[\partial_{\xi}^{(n-1)}\hat{\boldsymbol{A}}]\right]\Bigg)\\
=&\mathbb{T}_{\boldsymbol{Q}}[\partial_{\xi}^n\hat{\boldsymbol{A}}]+\sum_{k=1}^{n}\sum_{i=0}^{n-k}C_{i}^{(n-k)}\left[\partial_{\xi}^{(i)}\hat{\boldsymbol{\kappa}},\partial_{\xi}^{(n-k-i)}\mathbb{T}_{\boldsymbol{Q}}[\partial_{\xi}^{(k-1)}\hat{\boldsymbol{A}}]\right].\stepcounter{equation}\tag{\theequation}\label{myeq1}
\end{align*}
Noting that at $n=0$, the sum in equation \eqref{l42a} vanishes justifies the factor $(1-\delta_{n0})$. $\square$
\paragraph{Proposition 10:} Let $\hat{\boldsymbol{\kappa}}_\text{i}=\partial_{\xi}\boldsymbol{Q}_\text{i}.\boldsymbol{Q}_\text{i}^T$ and $\hat{\boldsymbol{\kappa}}_\text{f}=\partial_{\xi}\boldsymbol{Q}_\text{f}.\boldsymbol{Q}_\text{f}^T$ denote the curvature field corresponding to the initial and final configurations respectively. The updated curvature tensor and its derivatives are given by the recurrence-formula,
\begin{equation}\label{l43}
\begin{aligned}
\partial_{\xi}^n\hat{\boldsymbol{\kappa}}_\text{f}&=\partial_{\xi}^n\hat{\boldsymbol{\kappa}}_{+}+\mathbb{T}_{\boldsymbol{Q}_+}[\partial_{\xi}^n\hat{\boldsymbol{\kappa}}_\text{i}]+(1-\delta_{n0})\sum_{k=1}^{n}\sum_{i=0}^{n-k}C_{i}^{(n-k)}\left[\partial_{\xi}^{(i)}\hat{\boldsymbol{\kappa}}_{+},\partial_{\xi}^{(n-k-i)}\mathbb{T}_{\boldsymbol{Q}_+}[\partial_{\xi}^{(k-1)}\hat{\boldsymbol{\kappa}}_\text{i}]\right].
\end{aligned}
\end{equation}
\paragraph{Proof:} Using the Eq. \eqref{l33} and \eqref{l42}, we obtain updated curvature as:
\begin{equation}\label{l44}
\begin{aligned}
\hat{\boldsymbol{\kappa}}_\text{f}&=\partial_{\xi}(\boldsymbol{Q}_{+}.\boldsymbol{Q}_{i}).(\boldsymbol{Q}_{+}.\boldsymbol{Q}_{i})^T=\partial_{\xi}(\boldsymbol{Q}_{+}).\boldsymbol{Q}_{+}^T+\boldsymbol{Q}_{+}.\hat{\boldsymbol{\kappa}}_\text{i}.\boldsymbol{Q}_{+}^T=\hat{\boldsymbol{\kappa}}_{+}+\mathbb{T}_{\boldsymbol{Q}_+}[\hat{\boldsymbol{\kappa}}_\text{i}].
\end{aligned}
\end{equation}
Therefore,
\begin{equation}\label{l45}
\begin{aligned}
\partial_{\xi}^n\hat{\boldsymbol{\kappa}}_\text{f}&=\partial_{\xi}^n\hat{\boldsymbol{\kappa}}_{+}+\partial_{\xi}^n\mathbb{T}_{\boldsymbol{Q}_+}[\hat{\boldsymbol{\kappa}}_\text{i}].
\end{aligned}
\end{equation}
Substituting Eq. \eqref{l42a} obtained in \textit{Proposition 9} in the above result proves this proposition. $\square$

Appendix \ref{A63} presents a MATLAB code to update spatial curvature tensor and its derivatives in accordance with \textit{Proposition 10}.
\begin{figure}[htb!]\centering 
	\includegraphics[width=\textwidth]{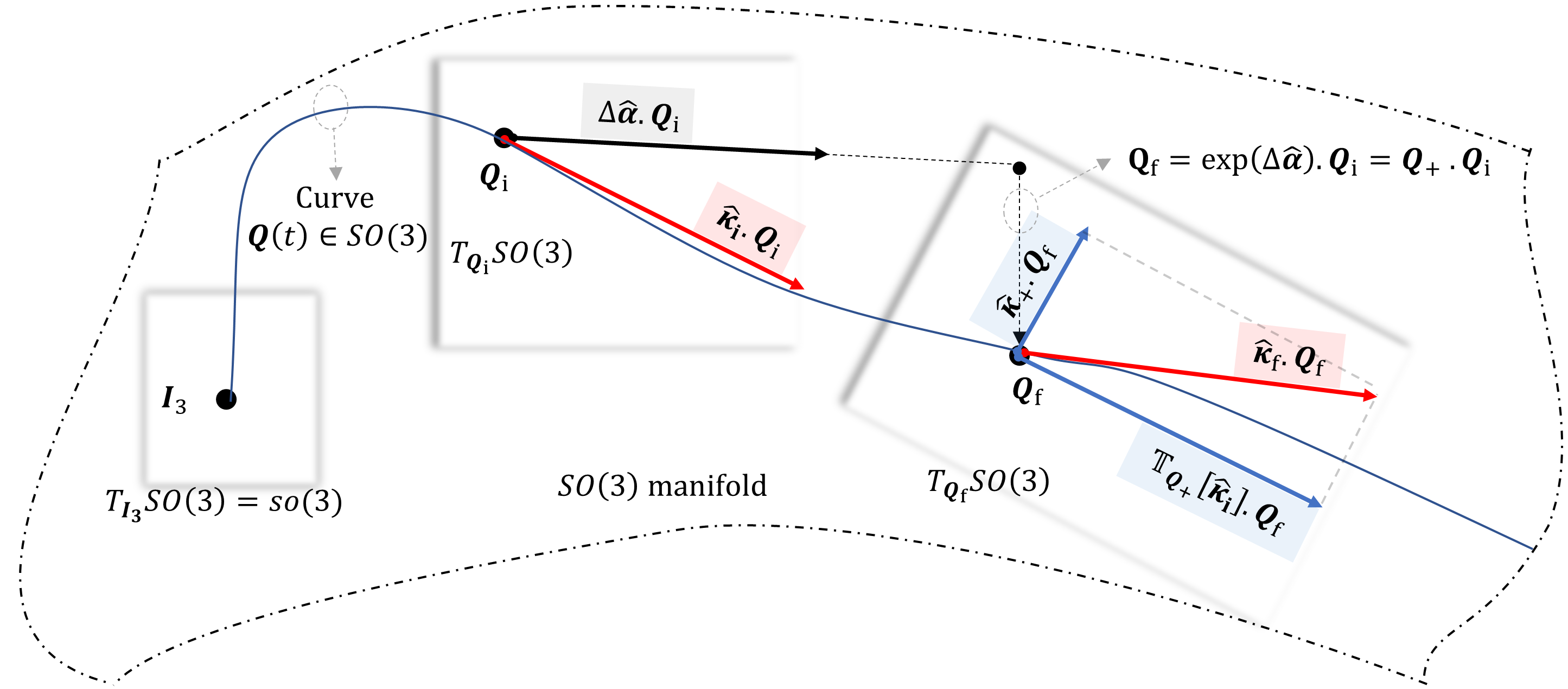}
	\caption{Geometric interpretation of curvature updating $\hat{\boldsymbol{\kappa}}_{\text{f}}=\hat{\boldsymbol{\kappa}}_{+}+\mathbb{T}_{\boldsymbol{Q}_+}[\hat{\boldsymbol{\kappa}}_\text{i}]$}
	\label{fig:Fig2}
\end{figure}

The curvatures $\hat{\boldsymbol\kappa}_{+}$ and their derivatives may be obtained using \textit{Proposition 7}. Once the spatial curvature and its derivatives are obtained, the derivative of material curvature and co-rotational derivative can be obtained using \textit{Proposition 1} or \textit{Corollary 3} or \textit{Corollary 4}. The reader should refer to figure \ref{fig:Fig2} for a geometric interpretation of curvature updating.
\FloatBarrier
\section{Conclusion}
\label{conclusion}
In this paper, we primarily focused our attention on the curvature tensor and its derivatives associated with any space curve framed by a general \textit{material frame}. Therefore, the results presented here are valid for any frame, including the Frenet-Serret and Bishop frames. This approach equivalently represents a curve on the $SO(3)$ manifold. In addition to discussing  the spatial and material forms of the curvature tensor, we have investigated the higher-order derivatives and co-rotation derivatives of these quantities. We have presented, for the first time to our knowledge, a closed-form formula for all higher-order derivative of the spatial curvature tensor. Finally, a time-updating algorithm for curvature (both spatial and material) and its derivatives (partial and co-rotational) was presented, which is particularly useful in practical problems like finite element formulation of geometrically-exact beams, among many other applications.

\paragraph{Acknowledgment:} Funding for this work was provided by the United States Army Corps of Engineers through the U.S. Army Engineer Research and Development Center Research Cooperative Agreement W912HZ-17-2-0024. We thank our colleague Professor J. S. Chen in Department of Structural Engineering at UC San Diego for providing us with
valuable advice and comments. 
\section{Appendix}
\label{appendix}
\subsection{MATLAB function to obtain spatial curvature and its derivatives}\label{A61}
We present a matlab function to obtain derivative of curvature upto order \hl{$\texttt{ord}$}. We assume that the rotation vector and its derivatives $\partial_{\xi}^n\boldsymbol{\uptheta}(\xi)$ for $0\leq n\leq (\hl{\texttt{ord}}+1)$ are known. We have $\partial_{\xi}^n\|\boldsymbol{\uptheta}\|=\partial_{\xi}^n\uptheta=\partial_{\xi}^n\sqrt{\boldsymbol{\uptheta}\cdot\boldsymbol{\uptheta}}$. The \textit{Gibbs vector} and its derivatives can then be obtained using $\partial_{\xi}^n\boldsymbol{\phi}(\xi)=\partial_{\xi}^n\left(\frac{\tan\frac{\uptheta}{2}}{\uptheta}{\boldsymbol\uptheta(\xi)}\right)$. Similarly, $\partial_{\xi}^n\overline{\phi}$ can be obtained using Eq. \eqref{l34}. 

Define the matrix \hl{\texttt{cur\_mat}} and \hl{\texttt{phi\_mat}} consisting of the curvature vector (to be obtained) and parametrizing vector $\boldsymbol{\phi}$ and its derivatives such that $m^{th}$ row gives $(m-1)$ derivative of the concerned quantity. Similarly, let \hl{\texttt{phi\_bar\_vec}} denote a vector consisting of the quantity $\overline{\phi}$ and its derivatives. 
\begin{equation*}
\begin{aligned}
&\hl{\texttt{cur\_mat}}=\left[\boldsymbol{\kappa};\partial_{\xi}\boldsymbol{\kappa};...;\partial_{\xi}^\texttt{ord}\boldsymbol{\kappa}\right];\\
&\hl{\texttt{phi\_mat}}=\left[\boldsymbol{\phi};\partial_{\xi}\boldsymbol{\phi};...;\partial_{\xi}^\texttt{ord}\boldsymbol{\phi};\partial_{\xi}^{\texttt{ord}+1}\boldsymbol{\phi}\right];\\
&\hl{\texttt{phi\_bar\_vec}}=\left[\overline{\phi},\partial_{\xi}\overline{\phi},...,\partial_{\xi}^{\texttt{ord}}\overline{\phi},\partial_{\xi}^{\texttt{ord}+1}\overline{\phi}\right].
\end{aligned}
\end{equation*}

Let \hl{\texttt{LB}} (for Lie Bracket), \hl{\texttt{hat}} and \hl{\texttt{unhat}} define MATLAB functions, with $\hl{\texttt{v1,v2,mat}}\in so(3)$ and $\hl{\texttt{v}}\in \mathbb{R}^3$ such that:
\begin{framed}
\begin{lstlisting}[language=Matlab]
function [mat] = LB(v1,v2)
mat=v1*v2-v2*v1;
\end{lstlisting}
\end{framed}
\begin{framed}
\begin{lstlisting}[language=Matlab]
function [mat] = hat(v)
mat=[0,-v(3),v(2);v(3),0,-v(1);-v(2),v(1),0];
\end{lstlisting}
\end{framed}
\begin{framed}
\begin{lstlisting}[language=Matlab]
function [v] = unhat(mat)
v=[mat(3,2),mat(1,3),mat(2,1)];
\end{lstlisting}
\end{framed}
The rotation tensor $\boldsymbol{Q}(\boldsymbol\uptheta)$, curvature and its derivatives can then be obtained using the function \hl{\texttt{cur\_der}} given \hl{\texttt{phi\_mat}} and \hl{\texttt{phi\_bar\_vec}} as the arguments.
\begin{framed}
\begin{lstlisting}[language=Matlab]
function [Q,cur_mat] = cur_der(phi_mat,phi_bar_vec)

%Obtain the highest order of curvature derivative required
ord=length(phi_bar_vec)-2;

%Obtain rotation tensor
Q=eye(3)+phi_bar_vec(1)*(hat(phi_matrix(1,:))+...
    hat(phi_matrix(1,:))*hat(phi_matrix(1,:)));

%Initialize kappa_mat 	
cur_mat=zeros(ord+1,3);

for m=0:ord
 for k=0:m
  integerTest=~mod((m-k+1)/2,1);
	
  if integerTest==1
	j_max=(m-k+1)/2-1;
  else
	j_max=floor((m-k+1)/2);
  end	
    dummy_mat=zeros(3,3);	
  for j=0:j_max
   b=((nchoosek(m-k,j)*((m-k-2*j+1)/(m-k-j+1)));
   dummy_mat=dummy_mat+...
   b*LB(hat(phi_mat(j+1,:)),hat(phi_mat(m-k-j+2,:)));
  end
	
  cur_mat(m+1,:)=cur_mat(m+1,:)+unhat(dummy_mat)+...
   nchoosek(m,k)*phi_bar_vec(k+1)*phi_mat(m-k+2,:);
 end
end
\end{lstlisting}
\end{framed}
\subsection{MATLAB function to obtain co-rotational derivatives, material curvature, its derivatives}\label{A62}
Let the co-rotational derivatives and material curvature and derivatives be defined as:
\begin{equation*}
\begin{aligned}
&\hl{\texttt{cur\_bar\_mat}}=\left[\overline{\boldsymbol{\kappa}};\partial_{\xi}\overline{\boldsymbol{\kappa}};...;\partial_{\xi}^\texttt{ord}\overline{\boldsymbol{\kappa}}\right];\\
&\hl{\texttt{cur\_tilde\_mat}}=\left[\tilde\partial_{\xi}{\boldsymbol{\kappa}};...;\tilde\partial_{\xi}^\texttt{ord}{\boldsymbol{\kappa}}\right];
\end{aligned}
\end{equation*}
Given the rotation tensor $\hl{\texttt{Q}}=\boldsymbol{Q}(\boldsymbol\uptheta)$, spatial curvature and its derivatives \hl{\texttt{cur\_mat}}, we can obtain \hl{\texttt{cur\_bar\_mat}} and \hl{\texttt{cur\_tilde\_mat}} using the function \hl{\texttt{cur\_bar\_tilde\_der}} with \hl{\texttt{Q}}; \hl{\texttt{cur\_mat}} and \hl{\texttt{ord}} as its arguments.
\begin{framed}
	\begin{lstlisting}[language=Matlab]
function [cur_bar_mat,cur_tilde_mat] 
                         = cur_bar_tilde_der(cur_mat,Q,ord)

%Step 1: Find der. of Q upto order ord using Prop 8.
%Q_der_mat(:,:,n) gives (n-1) der of Q. 
	
 Q_der_mat=zeros(3,3,ord+1);
 Q_der_mat(:,:,1)=Q;
	
 for n=1:ord
  for j=0:n-1
  Q_der_mat(:,:,n+1)=Q_der_mat(:,:,n+1)+...
   nchoosek(n-1,j)*hat(cur_mat(j+1,:))*Q_der_mat(:,:,n-j);
  end
 end
	
%Step 2: We use kappa_bar=transpose(Q).kappa and Corollary 4
	
 cur_bar_mat=zeros(3,ord);

 for n=1:ord
  for j=0:n        
  cur_bar_mat(:,n)=cur_bar_mat(:,n)+nchoosek(n,j)*...
  transpose(Q_der_mat(:,:,j+1))*cur_mat(n-j+1,:)';
  end
 end
 
 cur_bar=transpose(Q_der_mat(:,:,1))*cur_mat(1,:)';
 cur_bar_mat=[cur_bar,cur_bar_mat];
 cur_bar_mat=transpose(cur_bar_mat);
	
%Step 3: Evaluate co-rotational derivatives using Corollary 4

 cur_tilde_mat=zeros(ord,3);
 for n=1:ord+1
  cur_tilde_mat(n,:)=Q_der_mat(:,:,1)*cur_bar_mat(n+1,:)';
 end
	\end{lstlisting}
\end{framed}
\subsection{MATLAB function to update curvature and its derivatives}\label{A63}
Let the \textit{current rotation tensor} $\hl{\texttt{QI}}=\boldsymbol{Q}(\boldsymbol{\uptheta})=\boldsymbol{Q}_{i}$ (initial) be parametrized by rotation vector $\boldsymbol{\uptheta}(\xi)$. Let the \textit{current incremental rotation vector} (and its derivatives) be given as $\partial_{\xi}^n\Delta\boldsymbol\alpha$ for $0\leq n\leq (\hl{\texttt{ord}}+1)$, such that $\hl{\texttt{QP}}=\boldsymbol{Q}(\Delta\boldsymbol{\alpha})=\boldsymbol{Q}_{+}$. We reparametrize the two rotation vectors as $\boldsymbol{\phi}_\text{i}(\xi)=\frac{\tan\frac{\uptheta}{2}}{\uptheta}{\boldsymbol\uptheta(\xi)}$ and $\boldsymbol{\phi}_{+}(\xi)=\frac{\tan\frac{\|\Delta\boldsymbol\alpha\|}{2}}{\|\Delta\boldsymbol\alpha\|}{\Delta\boldsymbol\alpha(\xi)}$. Along the similar lines, we define the parameter $\overline\phi_\text{i}=2\cos^2\left(\frac{\uptheta}{2}\right)$ and $\overline\phi_\text{+}=2\cos^2\left(\frac{\|\Delta\boldsymbol\alpha\|}{2}\right)$.

\begin{equation*}
\begin{aligned}
&\hl{\texttt{phiI\_mat}}=\left[\boldsymbol{\phi}_{\text{i}};\partial_{\xi}\boldsymbol{\phi}_{\text{i}};...;\partial_{\xi}^\texttt{ord}\boldsymbol{\phi}_{\text{i}};\partial_{\xi}^{\texttt{ord}+1}\boldsymbol{\phi}_{\text{i}}\right];\\
&\hl{\texttt{phiI\_bar\_vec}}=\left[\overline{\phi}_{\text{i}},\partial_{\xi}\overline{\phi}_{\text{i}},...,\partial_{\xi}^{\texttt{ord}}\overline{\phi}_{\text{i}},\partial_{\xi}^{\texttt{ord}+1}\overline{\phi}_{\text{i}}\right];\\
&\hl{\texttt{phiP\_mat}}=\left[\boldsymbol{\phi}_{\text{+}};\partial_{\xi}\boldsymbol{\phi}_{\text{+}};...;\partial_{\xi}^\texttt{ord}\boldsymbol{\phi}_{\text{+}};\partial_{\xi}^{\texttt{ord}+1}\boldsymbol{\phi}_{\text{+}}\right];\\
&\hl{\texttt{phiP\_bar\_vec}}=\left[\overline{\phi}_{\text{+}},\partial_{\xi}\overline{\phi}_{\text{+}},...,\partial_{\xi}^{\texttt{ord}}\overline{\phi}_{\text{+}},\partial_{\xi}^{\texttt{ord}+1}\overline{\phi}_{\text{+}}\right].
\end{aligned}
\end{equation*} 

We assume that the above mentioned quantities are known. The corresponding rotation vectors, curvature and its derivatives can be obtained using the function defined in previous section as:
\begin{framed}
	\begin{lstlisting}[language=Matlab]
[QI,curI_mat] = cur_der(phiI_mat,phiI_bar_vec);
[QP,curP_mat] = cur_der(phiP_mat,phiP_bar_vec);
	\end{lstlisting}
\end{framed}
We define the function \hl{\texttt{TQ}} carrying out transport operation as defined by Eq. \eqref{l42}, such that for $\hl{\texttt{mat1,mat1}}\in so(3)$ and $\hl{\texttt{Q}}\in SO(3)$, we have
\begin{framed}
\begin{lstlisting}[language=Matlab]
function [mat2] = TQ(Q,mat1)
mat2=Q*mat1*transpose(Q);
\end{lstlisting}
\end{framed}

Let \hl{\texttt{QF}} represent the updated (final) rotation tensor and \hl{\texttt{kappaF\_mat}} the corresponding updated curvature and its derivatives that needs to be obtained. The function \hl{\texttt{cur\_updating}} defined yields obtain \hl{\texttt{QF}} and \hl{\texttt{kappaF\_mat}} given \hl{\texttt{QI,kappaI\_mat}} and \hl{\texttt{QP,kappaP\_mat}}. Before we present the function, we note from Eq. \eqref{l43} that to update, say upto $3^{\text{rd}}$ derivative of curvature vector (or tensor), we need to know the following quantities.
\begin{table*}[htb]
	\centering
	\begin{tabular}{|l|l|l|l|}
	    \hline \cellcolor{blue!10}$\mathbb{T}_{\boldsymbol{Q}_+}[\hat{\boldsymbol{\kappa}}_\text{i}]$& \cellcolor{blue!10}$\mathbb{T}_{\boldsymbol{Q}_+}[\partial_{\xi}\hat{\boldsymbol{\kappa}}_\text{i}]$ & \cellcolor{blue!10}$\mathbb{T}_{\boldsymbol{Q}_+}[\partial_{\xi}^2\hat{\boldsymbol{\kappa}}_\text{i}]$ & $\mathbb{T}_{\boldsymbol{Q}_+}[\partial_{\xi}^3\hat{\boldsymbol{\kappa}}_\text{i}]$ \\ \hline
		\cellcolor{blue!10}$\partial_{\xi}\mathbb{T}_{\boldsymbol{Q}_+}[\hat{\boldsymbol{\kappa}}_\text{i}]$& \cellcolor{blue!10}$\partial_{\xi}\mathbb{T}_{\boldsymbol{Q}_+}[\partial_{\xi}\hat{\boldsymbol{\kappa}}_\text{i}]$ & \cellcolor{blue!10}-& - \\ \hline
		\cellcolor{blue!10}$\partial_{\xi}^2\mathbb{T}_{\boldsymbol{Q}_+}[\hat{\boldsymbol{\kappa}}_\text{i}]$ & \cellcolor{blue!10}- & \cellcolor{blue!10}-& - \\ \hline
		\cellcolor{blue!10}$\partial_{\xi}^3\mathbb{T}_{\boldsymbol{Q}_+}[\hat{\boldsymbol{\kappa}}_\text{i}]$ & \cellcolor{blue!10}- & \cellcolor{blue!10}-& - \\\hline
	\end{tabular}
\end{table*}
The array \hl{\texttt{TQ\_mat}} stores these matrices for evaluation of the updated curvature. For $\hl{\texttt{ord}}=3$, the highlighted part of table below represents the array \hl{\texttt{TQ\_mat(:,:,i,j)}} where, $1\leq\texttt{i}\leq \hl{\texttt{ord}}+1$ and $1\leq\texttt{j}\leq \hl{\texttt{ord}}$.
\begin{framed}
\begin{lstlisting}[language=Matlab]
function [QF,curF_mat] = cur_updating(QI,QP,curP_mat,curI_mat)
ord=size(curI_mat,1)-1;

TQ_mat=zeros(3,3,ord+1,ord);

for n=1:ord+1

 if n==ord+1
  j_max=1;
 else
  j_max=ord+1-n;
 end
 
 for j=1:j_max
  
  if n==1
   TQ_mat(:,:,n,j)=TQ(QP,hat(curI_mat(n+j-1,:)));
  else
   dummy_mat=zeros(3,3);
   for k=1:n-1
    for i=0:(n-k-1)
    dummy_mat=dummy_mat+(nchoosek(n-k-1,i)
       *LB(hat(curP_mat(i+1,:)),TQ_mat(:,:,n-k-i,j+k-1)));
    end
   end 
   TQ_mat(:,:,n,j)=TQ(QP,hat(curI_mat(n+j-1,:)))+dummy_mat; 
  end 
   
 end   
end

%Update rotation tensor and curvature
QF=QP*QI;
curF_mat=zeros(4,3);
for n=1:4
  curF_mat(n,:)=curP_mat(n,:)+unhat(TQ_mat(:,:,n,1)) ;
end
\end{lstlisting}
\end{framed}
\newpage
\bibliographystyle{ieeetr}
\bibliography{reference}
\end{document}